\documentclass{article}
\usepackage[utf8]{inputenc}

\usepackage{etex}

\usepackage[linesnumbered,ruled,vlined,titlenotnumbered]{algorithm2e}
\SetCommentSty{mycommfont}
\SetKwInput{KwInput}{Input}         
\SetKwInput{KwOutput}{Output}

\usepackage[english]{babel}
\usepackage[utf8]{inputenc}
\usepackage[T1]{fontenc}

\usepackage{animate} %

\usepackage[most]{tcolorbox}

\usepackage{setspace}
\usepackage[absolute,overlay]{textpos}
\usepackage{tikz}
\usetikzlibrary{shapes.callouts}
\usetikzlibrary{decorations.text}
\newcommand\ind[1]{{1\kern-0.4em 1}_{\{#1\}}}

\usepackage{stmaryrd} %

\usetikzlibrary{matrix}

\usepackage{dsfont}
\usepackage{times}
\usepackage{amsmath}
\usepackage{amscd}
\usepackage{amsthm}
\usepackage[amsthm, amsmath, framed, thmmarks, thref]{}%
\usepackage{amssymb}
\usepackage{color}
\usepackage{enumerate}
\usepackage{todonotes}
\usepackage{mathtools}
\usepackage{theoremref}
\usepackage{empheq} %

\usepackage{mdwlist} %

\usepackage{multirow} %

\usepackage{graphicx}
\usepackage{caption}
\usepackage{subcaption}
\graphicspath{{FIGURES/}}

\usepackage{mathdots} %

\usepackage{amsfonts}

\usepackage{pifont}

\usepackage[all]{xy} %

\usepackage{cancel}

\usepackage{animate}

\usepackage{diagbox} %

\usepackage{multirow,array}

\usepackage{manfnt} %

\usepackage{hyperref}
\hypersetup{
    colorlinks=true,
    linkcolor=blue,
    filecolor=magenta, 
    urlcolor=cyan,
    citecolor=blue
    }

\usepackage{algorithmic}

\newcommand{\cA}{{\mathcal A}}
\newcommand{\cB}{{\mathcal B}}
\newcommand{\cC}{{\mathcal C}}

\newcommand{\cF}{{\mathcal F}}
\newcommand{\cG}{{\mathcal G}}

\newcommand{\cL}{{\mathcal L}}
\newcommand{\cK}{{\mathcal K}}
\newcommand{\cM}{{\mathcal M}}
\newcommand{\cN}{{\mathcal N}}

\newcommand{\cP}{{\mathcal P}}
\newcommand{\cQ}{{\mathcal Q}}

\newcommand{\cW}{{\mathcal W}}

\newcommand{\bbE}{{\mathbb E}}

\newcommand{\bbR}{{\mathbb R}}

\newcommand{\EE}{{\bbE}}

\newcommand{\RR}{{\bbR}}

\usepackage{xcolor}
\definecolor{mypurple}{RGB}{140,0,255}%
\definecolor{myred}{rgb}{255,0,0}
\definecolor{mydarkturquoise}{RGB}{0,206,209}
\definecolor{mydeeppink}{RGB}{255,20,147}
\definecolor{darkblue}{RGB}{0,0,140}
\definecolor{blue2}{RGB}{0,0,0}
\definecolor{middleblue}{RGB}{0,0,71}
\definecolor{light-gray}{gray}{0.9}
\definecolor{med-gray}{gray}{0.6}
\definecolor{ceruleanblue}{RGB}{51,102,187}
\definecolor{mediumblue}{RGB}{0,0,205}
\definecolor{darkorange}{RGB}{255,140,0}
\definecolor{orangered}{RGB}{255,69,0}
\definecolor{metallicorange}{RGB}{226, 99, 16}
\definecolor{dodgerblue}{RGB}{30, 144, 255}
\definecolor{seagreen}{RGB}{46, 139, 87}
\definecolor{med-gray}{gray}{0.6}
\definecolor{ProcessBlue}{cmyk}{1,0,0,0.40}
\definecolor{Black}{cmyk}{0,0,0,1}
\definecolor{Red}{cmyk}{0,1,1,0.2}
\definecolor{Green}{cmyk}{0.9,0,1,0}
\definecolor{Fuchsia}{cmyk}{0.47,0.91,0,0.06}
\definecolor{PineGreen}{cmyk}{0.92,0,0.59,0.30}
\definecolor{mydarkorange}{rgb}{1.0, 0.4, 0.0}

\definecolor{PurpleDeltaW}{RGB}{153,0,255}

\DeclareMathOperator{\diver}{div}
\DeclareMathOperator{\grad}{\nabla}

\DeclareMathOperator*{\argmax}{\arg\!\max}

\newcommand{\ctrl}{v}

\newtheorem{remark}{Remark}

\newtheorem{example}{Example}

\everymath{\textstyle}

\newtcolorbox{AlgoBox}[2][]{
                lower separated=false,
                colback=white!80!gray,
colframe=white!20!black,%
colbacktitle=white!30!gray,
coltitle=black,
enhanced,
attach boxed title to top left={xshift=0.5cm,
        yshift=-2mm},
title=#2,#1}
\newtcolorbox{AlgoBoxLighter}[2][]{
		bottom=0mm,
                lower separated=false,
                colback=white!90!gray,
colframe=white!20!black,%
colbacktitle=white!70!gray,
coltitle=black,
enhanced,
attach boxed title to top left={xshift=0.2cm,
        yshift=-2mm},
title=#2,#1}

\newcommand{\ctrldim}{k}
\newcommand{\dom}{\mathcal{Q}}
\newcommand{\domT}{\mathcal{Q}_{T}}

\newcommand{\forlongversion}[1]{#1}

\def\rd{\mathrm{d}}

\DeclarePairedDelimiterX{\inp}[2]{\langle}{\rangle}{#1, #2}

\newcommand\blfootnote[1]{%
  \begingroup
  \renewcommand\thefootnote{}\footnote{#1}%
  \addtocounter{footnote}{-1}%
  \endgroup
}

\usepackage{graphicx}
\usepackage{booktabs}
\usepackage{multirow,tabularx}
\newcommand{\tabincell}[2]{\begin{tabular}{@{}#1@{}}#2\end{tabular}}

\makeatletter
\def\blfootnote{\xdef\@thefnmark{}\@footnotetext}
\makeatother

\title{Deep Learning for Mean Field Optimal Transport
\blfootnote{{\bf Acknowledgements: } This project was realized during CEMRACS 2022. The authors would like to thank the CIRM for welcoming the CEMRACS 2022, the organizers of the CEMRACS 2022 for the opportunity to work on this project as well as their respective institutions. They are also grateful to the NYU-ECNU Institute of Mathematical Sciences at NYU Shanghai, the CIMPA fellowships program, and the ENS Rennes for their support. This work was supported in part through the NYUSH IT High Performance Computing resources, services, and staff expertise, as well as the ASCC Toubkal cluster resources.} 
}%
\author{
Sebastian Baudelet\thanks{Université Côte d'Azur, 28 Avenue de Valrose, 06103 Nice, France, sebastian.baudelet@univ-cotedazur.fr}, 
Brieuc Fr\'enais\thanks{IRMA UMR 7501, Université de Strasbourg, 7 Rue René Descartes, 67000 Strasbourg, France, brieuc.frenais@math.unistra.fr}, 
Mathieu Lauri\`ere\thanks{NYU-ECNU Institute of Mathematical Sciences at NYU Shanghai, 3663 Zhongshan Road North, Shanghai, 200062, China, mathieu.lauriere@nyu.edu}, 
Amal Machtalay\thanks{Mohammed VI Polytechnic University. Lot 660, Hay Moulay Rachid Ben Guerir, 43150, Morocco, amal.machtalay@um6p.ma}, 
Yuchen Zhu\thanks{Yale University, 433 Temple Street, New Haven, 06511, USA, yuchen.zhu@yale.edu}
}
\date{} 

\begin{document}

\maketitle

\begin{abstract} 
    Mean field control (MFC) problems have been introduced to study social optima in very large populations of strategic agents. The main idea is to consider an infinite population and to simplify the analysis by using a mean field approximation. These problems can also be viewed as optimal control problems for McKean-Vlasov dynamics. They have found applications in a wide range of fields, from economics and finance to social sciences and engineering. Usually, the goal for the agents is to minimize a total cost which consists in the integral of a running cost plus a terminal cost. In this work, we consider MFC problems in which there is no terminal cost but, instead, the terminal distribution is prescribed. We call such problems mean field optimal transport problems since they can be viewed as a generalization of classical optimal transport problems when mean field interactions occur in the dynamics or the running cost function. We propose three numerical methods based on neural networks. The first one is based on directly learning an optimal control. The second one amounts to solve a forward-backward PDE system characterizing the solution. The third one relies on a primal-dual approach. We illustrate these methods with numerical experiments conducted on two families of examples. 
\end{abstract}

\section{Introduction}

Mean field games (MFGs) have been introduced by Lasry and Lions~\cite{lasry2006jeuxi, lasry2006jeuxii, lasry2007mean} and Caines, Huang and Malham\'e~\cite{huang2006large, huang2007large} to approximate Nash equilibria in games with a very large number of players. At a high level, the main idea is to use a mean field approximation to represent the state of the population, and then to focus on the interactions between a single representative player and the distribution of the states of the other players. Mean field control (MFC)~\cite{bensoussan2013mean} relies on a similar approximation but aims at representing situations in which a large number of agents cooperate to minimize a common social cost. The problem can be interpreted as an optimal control problem for a McKean-Vlasov (MKV) stochastic differential equation (SDE) or an optimal control for a Kolmogorov-Fokker-Planck (KFP) partial differential equation (PDE). In the past decade, the analysis of both MFGs and MFC problems has been extensively developed, see e.g.~\cite{bensoussan2013mean} for an introduction to this topic, and ~\cite{carmona2018probabilistic} for a probabilistic viewpoint. 

In the most common setup, the players try to minimize a total cost which is composed of a running cost integrated over time and a terminal cost. These costs generally account for the efforts made to control the dynamics as well as the preferences for some states over others. Another class of models has been introduced, in which there is no terminal cost and instead the terminal distribution of the population is imposed as a constraint. Nash equilibria have been studied under the name of planning problem for mean field game. This class of problems has been analyzed mostly using PDE-based techniques~\cite{achdou2012mean,porretta2013planning, porretta2014planning,orrieri2019variational,graber2019planning,bertucci2021masterplanning}. In the special case of linear dynamics and a quadratic running cost in the control, the problem is related to the Schr\"odinger bridge problem, and equilibrium conditions can be phrased in terms of ordinary differential equations (ODEs)~\cite{chen2015optimal, chen2016relation, chen2015optimal2, chen2018optimal, chen2018steering}. 

This research direction is tightly connected to optimal transport (OT). Benamou and Brenier proposed in~\cite{benamou2000computational} a fluid mechanics framework for the $L^2$ Monge-Kantorovich mass transport problem. Many works built on this approach to relate optimal transport and optimal control problems for continuity equations. Of particular interest for MFGs is the work~\cite{cardaliaguet2013geodesics}, which clarified the link between geodesics for a class of distances between probability measures and a PDE system similar to the one arising in MFGs. For more background on OT, we refer the interested reader to the monographs~\cite{villani2009optimal, villani2021topics, santambrogio2015optimal,peyre2019computational,ambrosio2005gradient}. However, the solutions of MFGs correspond to Nash equilibria, and hence, in general, MFGs do not admit a variational structure. Furthermore, in many applications, it is not immediately clear to us why selfish players caring only about their individual costs would manage to agree and reach a target terminal distribution. Imposing a fixed terminal distribution seems more natural in the MFC setting, where the agents behave in a cooperative way to minimize the social cost. In the present work, we focus on such MFC with planning problems, in which a mean field of agents try to collectively minimize a social cost while ensuring that a fixed distribution is attained at the terminal time. 

Since the work of Benamou and Brenier~\cite{benamou2000computational}, several numerical methods have been investigated for similar problems, including MFGs with planning. Achdou et al. have proposed in~\cite{achdou2012mean} a method based on finite differences and Newton method to solve the PDE system of MFG with planning. Benamou and Carlier have used in~\cite{benamou2015augmented, benamou2016augmented} an Augmented Lagrangian method approach with the alternating direction method of multipliers to solve OT and MFG (without planning). Similar methods have been used in~\cite{andreev2017preconditioning,achdou2016meancongestion2} to solve MFGs and MFC problems (still without planning). Benamou et al. proposed in~\cite{benamou2019entropy} a method to solve MFG with planning through entropy regularization and Sinkhorn algorithm. 

Recently, several deep learning methods have been proposed to solve high-dimensional optimal control problems and PDEs, such as the DeepBSDE method~\cite{han2018solving, han2017deep, han2022learning}, the Deep Galerkin Method~\cite{sirignano2018dgm} and physics-informed neural networks~\cite{raissi2019physics}. Some of these methods have been extended to MFGs and MFC problems. In particular, \cite{alaradi2018solving,carmona2019convergence} proposed deep learning methods to solve the PDE systems arising in mean field problems, \cite{carmona2021convergence,germain2022numerical,fouque2020deep} introduced deep learning methods for differential MFC problems. Ruthotto et al. introduced a deep learning method for variational MFG with degenerate diffusion in~\cite{ruthotto2020machine}. Lin et al. introduced a deep learning method in ~\cite{lin2020apac} that utilizes the primal-dual relationship of variational MFG. Cao et al. noticed a connection between MFGs, generative adversarial networks and OT in~\cite{cao2020connecting}. We refer the interested reader to e.g.~\cite{carmona2021deep,germain2021neural,hu2022recent} for recent surveys on this topic. The work most related to ours is the work of Liu et al. in ~\cite{liu2022deep}, where they considered the planning problems in a class of MFGs based on a generalized version of the Schr\"odinger bridge problem and proposed a neural network-based numerical method to solved it.

The main goal of this paper is to propose numerical methods based on deep learning to solve MFC problems with planning constraint, that we will call mean field optimal transport problems. To the best of our knowledge, the theory remains to be investigated in detail, and this is beyond the scope of the present work. Here, we proceed formally when needed, and we focus on the numerical aspects using machine learning tools. The rest of the paper is organized as follows. In Section~\ref{sec:defpb}, we introduce the problem and discuss several examples. In Section~\ref{sec:methods}, we describe three numerical methods, each based on a different approach for the problem. In Section~\ref{sec:num-exp}, we present numerical results on several benchmark problems.

\section{Definition of the problem}
\label{sec:defpb}

Before presenting the mean field optimal transport problem, let us first recall the definition of a typical mean field control problem. Let $T$ be a time horizon. Let $\dom = \RR^d$ and $\domT = [0,T] \times \dom$ denote the space domain and the time-space domain. Denote by $\cP_2(\dom)$ the set of square-integrable probability measures on $\dom$. Let $f: \dom \times \cP_2(\dom) \times \RR^k \to \RR$ be a running cost function, $g: \dom \times \cP_2(\dom) \to \RR$ be a terminal cost function, $b: \dom \times \cP_2(\dom) \times \RR^k \to \RR^d$ be a drift function and $\sigma \in \RR$ be a non-negative constant diffusion coefficient.  
In a classical MFC problem with given initial distribution $\rho_0$ in $\cP_2(\dom)$, the goal is to find a feedback control $\ctrl^*: \domT\to \RR^\ctrldim$ minimizing:
\begin{align}
\label{eq:def-J-MFC}
	J^{MFC}: \ctrl \mapsto \EE \left[\int_0^T f(X_t^{\ctrl}, \mu^{\ctrl}(t), \ctrl(t,X_t^{\ctrl}) ) \mathrm{d}t + g(X_T^{\ctrl}, \mu^{\ctrl}(T)) \right]
\end{align}
where $\mu^{\ctrl}(t)$ is the distribution of $X_t^{\ctrl}$, under the constraint that the process $X^{\ctrl} = (X_t^{\ctrl})_{t \ge 0}$ solves the SDE
\begin{align}
\label{eq:dyn-X-general-MFC}
\begin{cases}
    X_0^{\ctrl} \sim \rho_0
    \\
	\mathrm{d} X_t^{\ctrl} = b(X_t^{\ctrl}, \mu^{\ctrl}(t), \ctrl(t, X_t^{\ctrl})) \mathrm{d}t + \sigma \mathrm{d} W_t, \qquad t \ge 0,
\end{cases}
\end{align}
where $W$ is a standard $d$-dimensional Brownian motion. It would also be interesting to consider open-loop controls, but since we are motivated by numerical applications, we restrict our attention to feedback controls. The cost~\eqref{eq:def-J-MFC} can be interpreted either as the expected cost for a single representative player, or as the average cost for the whole population, which we refer to as the social cost.

In this work, we are interested in a modified version of the above problem, where instead of having a terminal cost, a terminal distribution is imposed. This type of problem encompasses optimal transport as a special case, but it may incorporate mean field interactions in the drift and the running cost. For this reason, we will refer to this class of problems as mean field optimal transport (MFOT for short).\footnote{By analogy with MFG of planning type, we could also call such problems ``MFC of planning type''. But referring to ``optimal transport'' seems clearer so we will stick to the MFOT terminology.} 
Given two distributions $\rho_0$ and $\rho_T \in \cP_2(\dom)$, the goal is to find a feedback control $\ctrl^*: \domT\to \RR^\ctrldim$ minimizing
\begin{align}
\label{eq:def-J-MFOT}
    J^{MFOT}: \ctrl \mapsto \EE \left[\int_0^T f(X_t^{\ctrl}, \mu^{\ctrl}(t), \ctrl(t,X_t^{\ctrl}) ) \mathrm{d}t \right],
\end{align}
where $\mu^{\ctrl}(t)$ is the distribution of $X_t^{\ctrl}$, under the constraint that the process $X^{\ctrl} = (X_t^{\ctrl})_{t \ge 0}$ solves the SDE
\begin{align}
\begin{cases}
\label{eq:dyn-MFOT}
    X_0^{\ctrl} \sim \rho_0, \qquad X_T^{\ctrl} \sim \rho_T
    \\
    d X_t^{\ctrl} = b(X_t^{\ctrl}, \mu^{\ctrl}(t), \ctrl(t, X_t^{\ctrl})) \mathrm{d}t + \sigma \mathrm{d} W_t, \qquad t \ge 0.
\end{cases}
\end{align}
We stress that the terminal constraint implicitly restricts the class of admissible controls since we are interested in minimizing only over controls $\ctrl$ that make $X_T^{\ctrl}$ have distribution $\rho_T$.

We now present a few useful examples, some of which will be revisited in the numerical experiments (see Section~\ref{sec:num-exp}).

\begin{example}[Optimal transport]
When $b(x,\mu,a) = a$, $f(x,\mu,a) = \frac{1}{2} a^\top a$ and $\sigma=0$, the MFOT problem reduces to a \emph{standard OT} problem. See e.g.~\cite{benamou2000computational}.
\end{example}

\begin{example}[Linear-quadratic]
\label{ex:lq-example-description}
    Take $b(x,\mu,a) = Ax + \bar A \bar \mu + B a$, $f(x,\mu,a) = x^\top Q x + \bar\mu^\top \bar Q \bar \mu + a^\top R a$, and $g(x,\mu) = x^\top Q_T x + \bar\mu^\top \bar Q_T \bar \mu$, where $\bar\mu = \int \xi \mu(\mathrm{d}\xi)$, where $A, \bar A, B, Q, \bar Q, R, Q_T$ and $\bar Q_T$ are matrices of suitable sizes. In this setting, the MFC problem has an explicit solution, up to solving a forward-backward system of ODEs.  Furthermore, if the initial distribution is Gaussian, then the optimal flow of distribution remains Gaussian. See e.g.~\cite[Chapter 6]{bensoussan2013mean}. To the best of our knowledge, in the MFOT setting, a similar result is available in the literature only when $Q = \bar Q = 0$, which corresponds to the Schr\"odinger bridge problem. See~\cite[Section 7.1]{chen2018steering}.
\end{example}

\begin{example}[Crowd motion with congestion]
\label{ex:congestion}
\label{ex:congestion-example-description}
    Take $b(x,\mu,a) = a$, $f(x,\mu,a) = (c+\rho \star \mu(x))^\gamma |a|^2 + \ell(x, \mu(x))$, where $c \ge 0$ is a constant, $\rho$ is a regularizing kernel and $\star$ denotes the convolution. For $\gamma=0$, the model is linear-quadratic in the control. If $\gamma>0$, the cost of moving increases with the density surrounding the agent, which represents the fact that the ``energy'' spent to move is higher in regions with higher density. This models a \emph{congestion effect}. The last term in $f$ can be used to represent \emph{crowd aversion} if $\ell$ is increasing with respect to $\mu(x)$, and it can be used to represent spatial preferences by taking for instance $\ell(x, \mu(x)) = |x_* - x|^2$, where $x_*$ is a preferred position. The terminal cost $g$ can also be used to represent crowd aversion or spatial preferences. See e.g.~\cite{achdou2015system,achdou2016meancongestion1} for more details on the analysis of the MFC PDE system for this class of models and~\cite{achdou2016meancongestion2} for numerical aspects. When $\ell=0$ and $\sigma=0$, the corresponding MFOT problem has been studied e.g. in~\cite{cardaliaguet2013geodesics}. Similar models have also been studied in the context of MFGs, see e.g.~\cite{achdouporretta2018mean,achdoulasry2019mean}.
\end{example}

\section{Numerical methods} 
\label{sec:methods}
In this section, we introduce three different numerical methods to solve MFOT. Section~\ref{sec:method1} introduces a direct approach to solve a MFC problem that approximates the MFOT problem. Section~\ref{subsec:dgm} discusses the Deep Galerkin Method (DGM) to solve the underlying PDE system that characterizes the optimal solution to MFOT, which is composed of a coupled Hamilton-Jacobi-Bellman equation and a Kolomogrov-Fokker-Planck equation. Section~\ref{subsec:deepadmm} introduces the DeepADMM algorithm that solves a variational reformulation of the MFOT problem based on an augmented Lagrangian approach.
\subsection{Direct approach for the optimal control formulation}
\label{sec:method1}

We first introduce the direct approach, which does not require any derivation of optimality conditions. In order to make the problem numerically tractable, we make approximations on several levels. Motivated by the deep learning method for MFC problems proposed in~\cite{carmona2019convergence} (see also the first algorithm in~\cite{carmona2021deep}), we first approximate the MFOT problem~\eqref{eq:def-J-MFOT} by an MFC problem in which a terminal penalty is incurred based on the distance between the terminal distribution and the target distribution. We can then apply the algorithm of~\cite{carmona2019convergence}, which trains a neural network to learn the optimal control of the MFC problem. This method itself relies on three approximations. 

\subsubsection{Problem Approximation}

Instead of directly tackling the MFOT problem~\eqref{eq:def-J-MFOT}, we first consider the following MFC problem as an approximation of the original problem: Find a feedback control $\ctrl^*: \domT\to \RR^\ctrldim$ minimizing~\eqref{eq:def-J-MFC} under the constraint~\eqref{eq:dyn-X-general-MFC} when the terminal cost is:
\begin{equation}
    \label{eq:M1-penalty}
    g(x, \mu) = G(\mathcal{W}_2(\mu, \rho_T)), \qquad \mu \in \cP_2(\cQ).
\end{equation}
where $G:\RR_+ \to \RR_+$ is an increasing function and $\mathcal{W}_2$ denotes the Wasserstein distance on $\cP_2(\dom)$. A typical example that we will use in the experiments is a linear function. The purpose of introducing $G(\mathcal{W}_2(\mu, \rho_T))$ is to add a penalty that enforces the planning constraint for the terminal distribution. Here we focus on the Wasserstein distance because of its connection with optimal transport, see e.g.~\cite{benamou2000computational,santambrogio2015optimal}, although other similarity measures could be used. In our numerical experiments, we will take an increasing linear function for $G$. %

Then, we use the following approximations:
\begin{itemize}
    \item Since it is not possible to optimize overall feedback controls, we restrict the space of controls to the space of neural networks with a given architecture. We will denote by $\ctrl_\theta$ a representative neural network of this class with parameter $\theta$. The problem becomes a finite-dimensional optimization problem, in which the goal is to find a value for the parameter $\theta$ that minimizes the loss $J(\ctrl_\theta)$, i.e., the total cost of the MFC problem when using control $\ctrl_\theta$. 
    \item Since it is not possible to represent the mean field state $\mu^{v_\theta}(t)$ or to compute its evolution exactly, we approximate it by the empirical distribution $\bar\mu^{N,v_\theta}(t) = \frac{1}{N} \sum_{i=1}^N \delta_{X^{i,\ctrl_\theta}_t}$, where each $X^{i,v_\theta}$ is a solution of,
    \begin{align}
    \label{eq:dyn-X-general-MFC-N}
    \begin{cases}
        X_0^{i,\ctrl_\theta} \sim \rho_0 \quad \hbox{ i.i.d.}
        \\
    	d X_t^{i,\ctrl_\theta} = b(X_t^{i,\ctrl_\theta}, \bar\mu^{N,\ctrl_\theta}(t), \ctrl_\theta(t, X_t^{i,\ctrl_\theta})) \mathrm{d}t + \sigma \mathrm{d} W_t^i, \qquad t \ge 0,
    \end{cases}
    \end{align}
    where $(W^i)_{i=1,\dots,N}$ is a family of $N$ independent $d$-dimensional Brownian motions, which represent idiosyncratic noises affecting each particle independently. All the SDEs are based on the same control function $\ctrl_\theta$. 
    \item Last, in order to be able to compute these dynamics using Monte Carlo simulations, we discretize the time variable $t$. Letting $N_T$ be a number of regular time steps of length $\Delta t = T / N_T$, we replace the interval $[0,T]$ by the time steps $\{t_0=0, t_1=\Delta t, \dots, t_{N_T} = N_T \Delta t \}$. The time steps are $t_n = n \Delta t$, $n=0,\dots,N_T$. We then approximate the SDE system~\eqref{eq:dyn-X-general-MFC-N} using an Euler-Maruyama scheme. The
    family of trajectories $((X^{i,\ctrl_\theta}_t)_{t \in [0,T]})_{i=1,\dots,N}$ is approximated by the family of sequences $((X^{i,\ctrl_\theta,N_T}_{t_n})_{n=0,\dots,N_T})_{i=1,\dots,N}$ satisfying:
    \begin{equation}\label{eq:dyn-X-general-MFC-N-Deltat}
    \begin{cases}
        X^{i,\ctrl_\theta,N_T}_0 \sim \rho_0 \quad \hbox{ i.i.d.}
        \\
    	X^{i,\ctrl_\theta,N_T}_{t_{n+1}}=X^{i,\ctrl_\theta,N_T}_{t_n}+b(X^{i,\ctrl_\theta,N_T}_{t_n},\bar\mu^{N,v_\theta,N_T}_{t_n}, \ctrl_\theta(t_n,X^{i,\ctrl_\theta,N_T}_{t_n}))\Delta t + \sigma \Delta W^i_n,
    \end{cases}
    \end{equation}
    where $\bar\mu^{N,v_\theta,N_T}_{t_n} = \frac{1}{N} \sum_{i=1}^N \delta_{X^{i,\ctrl_\theta,N_T}_{t_n}}$ is the empirical distribution associated with the samples $X^{i,\theta,N_T}_{t_n}$. Here, $(\Delta W^i_n)_{i=1,\dots,N,n=0,\dots,N_T-1}$ are independent Gaussian random variables with variance $\Delta t$. 
\end{itemize}
To summarize, the new problem is to find $\theta^*$ minimizing:
\begin{equation*}
    J^{N,N_T}(\theta) 
    = \mathbb{E}\left[ \dfrac{1}{N}\sum_{i=1}^N\sum_{n=0}^{N_T-1}f(X^{i,\theta,N_T}_{t_n}, \bar\mu^{N,\theta,N_T}_{t_n},\ctrl_\theta(t_n,X^{i,\theta,N_T}_{t_n}))\Delta t+g(X^{N,\theta,N_T}_{T},\bar\mu^{N,\theta,N_T}_{T}) \right]
\end{equation*}
subject to the dynamics~\eqref{eq:dyn-X-general-MFC-N-Deltat}. The full analysis of this problem and its rigorous connection with the original MFOT problem~\eqref{eq:def-J-MFOT} is beyond the scope of this paper and is left for future work. We expect the control $v_{\theta^*}$, with the parameter value that is optimal for the above problem, to be approximately optimal for~\eqref{eq:def-J-MFOT}, under suitable assumptions on $b$ and $f$. In particular, $b$ and $f$ should probably depend smoothly on the distribution so that they can be evaluated in a meaningful way at the empirical distribution $\bar\mu^{N,\theta,N_T}_{t_n}$.

\subsubsection{Description of the algorithm}
\label{subsec:method1-algo}

\textbf{Optimization method. } To find an approximate minimizer, we use stochastic gradient descent (SGD) or one of its variants. At iteration $k$, we have a parameter $\theta_k$ that we wish to update. We sample the initial positions $(X^{i,\ctrl_{\theta_k},N_T}_0)_{i=1,\dots,N}$ and the Brownian motion increments $(\Delta W^i_n)_{i=1,\dots,N,n=0,\dots,N_T-1}$. We then compute the empirical cost for this realization of the $N$-particle population, and use its gradient with respect to $\theta_k$ to update the parameter. In other words, we apply SGD to the following loss function:
$$
    \mathcal{L}(\theta) 
    = J^{N,N_T}(\theta) 
    = \mathbb{E}_{S}[\mathcal{L}(\theta; S)],
$$
with:
$$
    \mathcal{L}(\theta; S) 
    =  \dfrac{1}{N}\sum_{i=1}^N\sum_{n=0}^{N_T-1}f(X^{i,\theta,N_T}_{t_n}, \bar\mu^{N,\theta,N_T}_{t_n},\ctrl_\theta(t_n,X^{i,\theta,N_T}_{t_n}))\Delta t+g(X^{N,\theta,N_T}_{T},\bar\mu^{N,\theta,N_T}_{T})
$$
where $S = \left((X^{i,\ctrl_\theta,N_T}_0)_{i=1,\dots,N}, (\Delta W^i_n)_{i=1,\dots,N,n=0,\dots,N_T-1} \right)$ denotes one random sample.

\textbf{Computation of the Wasserstein distance. } 
As shown in~\eqref{eq:M1-penalty}, the new problem we considered involves a Wasserstein distance between two continuous distributions, namely, the mean field distribution at terminal time $\mu^\ctrl_T$ and the target distribution $\rho_T$.  This is in general hard to compute. However, in our implementation, the mean field distribution is approximated by an empirical distribution obtained by Monte Carlo simulations, as is explained above. We then sample the same number of points from the target distribution and compute the Wasserstein distance between the two empirical distributions. This is done in the following way. 
Let $X$ and $Y$ be two sets of $N$ points each sampled from distributions $\mu, \nu$, $M_p$ the distance matrix, $(M_p)_{ij} = |X_i - Y_j|^p $, and the following set:
\begin{align}
\label{eq:set-U-Wasserstein}
        U_N = \left\lbrace A \in \mathbb{R}^{N \times N} \Big| \sum\limits_{j = 1}^N A_{ij} = \sum\limits_{i = 1}^{N} A_{ij} = \dfrac{1}{N}\right\rbrace.
\end{align}
Then
\begin{align*}
        \Bigl(\mathcal{W}_p\left(\mu,\nu \right)\Bigr)^{p} = \lim\limits_{N \to \infty} \min\limits_{T \in U_N} \left\langle T,M_p \right\rangle.
\end{align*}

In order to efficiently compute the Wasserstein distance, we follow the algorithm proposed by Cuturi in~\cite{cuturi2013sinkhorn}. We consider an extra entropy regularization of the following form. Let $\alpha > 0$, we want to find $T^*_\alpha$,  which is the solution to the following program:
\begin{align*}
    \min\limits_{T \in U_N} \left\langle T,M_p \right\rangle - \alpha \left\langle T \log(T), 1 \right\rangle.
\end{align*}
Optimality conditions and Sinkhorn-Knopp \cite{sinkhorn1967concerning} theorem give us the existence and uniqueness of the solution, as well as a unique decomposition of $T_{\alpha}^{*}$ using two vectors $u$ and $v$ such that:
\begin{align*}
    T^*_\alpha = \operatorname{Diag}(u)\exp\left(-\dfrac{M_p}{\alpha}\right)\operatorname{Diag}(v).
\end{align*}
We can then compute $u$ and $v$ with Sinkhorn-Knopp algorithm. Further explanations on this algorithm can be found in \cite{cuturi2013sinkhorn}. This method allows for fast computations and is easy to export to greater dimensions, at the cost of adding a layer of approximation due to the extra parameter $\alpha$. It can be noticed that as $\alpha$ tends to zero, the regularized solution tends to the solution of discrete optimal transport. In practice, reducing $\alpha$ to zero increases the number of iterations required for Sinkhorn algorithm to converge. However, in our numerical experiments we usually obtain good results with a small but non-zero $\alpha$.  %

\begin{remark}
    Notice that using our approach, we have one empirical distribution and one continuous distribution. Indeed, we have the empirical distribution obtained by Monte Carlo simulation and the target distribution $\rho_T$, which is generally given by a closed-form formula for its density. We could thus try to use the designated methods, such as Semi-discrete Optimal transport~\cite{merigot2021optimal}. While being more accurate, these methods do not scale well in higher dimensions compared to Sinkhorn's alternative.
\end{remark}

\textbf{Terminal penalty. } 
In our implementation, we take $G$ as a linear function $G(r) = C_{W} r$, where $C_W$ is a positive constant that weighs the importance of the terminal penalty in comparison with the running cost. This leads to a trade-off between minimizing the running cost and satisfying the terminal constraint. We noticed that when $C_W$ is too small, the algorithm minimizes the running cost without much consideration for the terminal condition and hence the terminal distribution is far from the target distribution. Therefore, the penalization has to be a significant component of the total loss if we want the terminal planning constraint to be approximately satisfied with good accuracy.

\subsection{Deep Galerkin Method for the PDE system} 
\label{subsec:dgm}

We now turn our attention to a method based on solving a forward-backward PDE system that characterizes the solution. We first discuss the PDE system and then use a deep learning method to solve this system.

\subsubsection{PDE system for MFOT}
As recalled above, in a standard MFC, the whole population uses a given feedback control $v$. Assuming that the distribution $\mu^v_t = \cL(X^v_t)$ of a representative agent with dynamics~\eqref{eq:dyn-X-general-MFC} admits a smooth enough density $m_t$, the latter satisfies the Kolmogorov-Fokker-Planck (KFP) PDE: 
\begin{equation}
\label{eq:KFPPDE-v}
    \begin{cases}
        &\frac{\partial m}{\partial t}(t,x) - \nu \Delta m(t,x) + \diver \bigl( m(t, x) b(x,m(t,\cdot),v(t, x)) \bigr) = 0 \qquad t \in (0,T], x \in \dom \nonumber \\
    & m(0,x) = m_{0}(x), \qquad x \in \dom,
    \end{cases}
\end{equation}
where $m_0$ is the density of the initial distribution $\rho_0$ and $\nu = \frac{\sigma^2}{2}$. The MFC problem~\eqref{eq:def-J-MFC} can then be viewed as an optimal control problem driven by the above KFP PDE. Under suitable conditions, the optimal control can be characterized through an adjoint PDE, which can be derived for instance via calculus of variations. See \textit{e.g.}~\cite[Chapter 4]{bensoussan2013mean} for more details.

Let $H: \dom \times L^2(\dom) \times \RR^d  \to  \RR$ be the Hamiltonian of the control problem faced by an infinitesimal agent in the first point above, which is defined by:
\begin{align}
\label{eq:def-hamiltonian}
    H: (x,m,p) \mapsto H(x,m,p) = \max_{\ctrl \in \RR^\ctrldim} \{ -L(x,m,\ctrl,p) \},
\end{align}

where $m$ denotes the density of $\mu$, $L: \dom \times L^2(\dom) \times \RR^\ctrldim \times \RR^d \to  \RR$ is the Lagrangian, defined by:
\begin{align}
\label{eq:def-lagrangian}
    L: (x,m,\ctrl,p) \mapsto L(x,m,\ctrl,p) =
    f(x,m,\ctrl) + \langle b(x,m,\ctrl) , p \rangle.
\end{align}

A necessary condition for the existence of an optimal control $\ctrl^* $ is that:
$$
  \ctrl^*(t,x) = \argmax_{\ctrl \in \RR^\ctrldim} \big\{ -L(x, m(t,\cdot), \ctrl, \nabla u(t,x))  \big\},
$$
where $(u,m)$ solve the
following system of partial differential equations: %
\begin{subequations} 
     \begin{empheq}[left=\empheqlbrace]{alignat=2}
	0 &=
	\displaystyle
	-\frac{\partial u} {\partial t} (t,x) - \nu \Delta u(t,x) + H( x, m(t,\cdot), \nabla u(t,x)) 
	\notag
	\\
	&\qquad + \int_{\dom} \frac{\partial H} {\partial m}(\zeta, m(t,\cdot), \nabla u(t, \zeta))(x) m(t,\zeta) \mathrm{d}\zeta  , 
	&& \hbox{ in } (0,T] \times \dom,
\\
	0 &= 
	\displaystyle \frac{\partial m} {\partial t} (t,x)  - \nu \Delta  m(t,x)
	- \diver\Bigl( m(t,x) \partial_p H(x, m(t, \cdot),\nabla u(t,x))\Bigr) ,
	&& \hbox{ in } [0,T) \times \dom,
	\\
  	& u(T,x) = g (x, m(T,\cdot)) 
	+ \int_{\dom}
\frac{\partial g} {\partial m} (\zeta, m(T,\cdot))(x) m(T,\zeta) \mathrm{d}\zeta, 
	&& \hbox{ in }  \dom,
	\\
	& m(0,x) = m_0(x),  && \hbox{ in } \dom.
\end{empheq}
\end{subequations}
The partial derivatives with respect to $m$ appear in the backward PDE due to the fact that the population distribution changes when the control changes.  These partial derivatives with respect to $m$ should be understood in the following sense: if $\varphi: L^2(\RR^d) \to \RR$ is differentiable, 
$$
	\frac{\mathrm{d}}{\mathrm{d}\varepsilon} \varphi(m + \varepsilon \tilde m)(x)_{\big| \varepsilon = 0} = \int_{\RR^d} \frac{\partial \varphi}{\partial m}(m)(\zeta) \tilde m(\zeta) \mathrm{d}\zeta.
$$

 We refer to e.g.~\cite[Chapter 4]{bensoussan2013mean} for more details and for the derivation using calculus of variations, which clarifies why the partial derivatives with respect to $m$ appear. If the cost functions and the drift function depend on the density only locally (i.e., only on the density at the current position of the agent), $\frac{\partial}{\partial m}$ becomes a derivative in the usual sense. 

  In this PDE system, $m$ plays the role of the MFC problem's state. The forward equation is a Kolmogorov-Fokker-Planck (KFP) equation which describes the evolution of the mean field distribution.  The other unknown function, $u$, plays the role of an adjoint state. Although the backward PDE has the form of a Hamilton-Jacobi-Bellman (HJB) equation, $u$ cannot, in general, be interpreted as the value function associated to problem~\eqref{eq:def-J-MFC} because the value function depends on the population distribution, see e.g.~\cite{lauriere2014dynamic,bensoussan2015master}. We refer the interested reader to e.g.~\cite[Chapters 3 and 4]{bensoussan2013mean} for the comparison with the MFG PDE system, in which the terms involving a derivative with respect to $m$ are absent, and $u$ can be interpreted as the value function of an infinitesimal player. 
 
Now, for the MFOT problem, we can proceed formally in a similar way. We derive an analogous PDE system, except that the terminal condition for $u$ disappears, and a terminal condition for $m$ is added to the system. More precisely, we (formally) obtain the following PDE system: 
\begin{subequations} 
\label{eq:MFOT-pde-system}
     \begin{empheq}[left=\empheqlbrace]{alignat=2}
	0 &=
	\displaystyle
	-\frac{\partial u} {\partial t} (t,x) - \nu \Delta u(t,x) + H( x, m(t,\cdot), \nabla u(t,x)) 
	\notag
	\\
	&\qquad + \int_{\dom} \frac{\partial H} {\partial m}(\zeta, m(t,\cdot), \nabla u(t, \zeta))(x) m(t,\zeta) \mathrm{d}\zeta  , 
	&& \hbox{ in } (0,T] \times \dom,
\\
	0 &= 
	\displaystyle \frac{\partial m} {\partial t} (t,x)  - \nu \Delta  m(t,x)
	- \diver\Bigl( m(t,x) \partial_p H(x, m(t,\cdot),\nabla u(t,x))\Bigr) ,
	&& \hbox{ in } [0,T) \times \dom,
	\\
	& m(0,x) = m_0(x),  \qquad m(T,x) = m_T(x) && \hbox{ in } \dom.
\end{empheq}
\end{subequations}
where $m_0$ and $m_T$ are respectively the densities of $\rho_0$ and $\rho_T$. To the best of our knowledge, this PDE system has not been derived nor analyzed in a general setting. Notice that even the existence of a solution is a non trivial question due to the fact that there is both an initial and a terminal constraint on the density. However, this system has been analyzed in special cases corresponding to optimal transport~\cite{benamou2000computational,cardaliaguet2013geodesics} or to MFGs with planning~\cite{achdou2012mean,porretta2013planning,porretta2014planning,graber2019planning}. In the numerical examples of Section~\ref{sec:num-exp}, we will mostly focus on cases that have been previously studied, such as standard optimal transport or MFOT with congestion effects captured by the running cost.

\subsubsection{Description of the algorithm}

To solve the PDE system~\eqref{eq:MFOT-pde-system}, we follow the idea of the Deep Galerkin Method (DGM) introduced by Sigignano and Spiliopoulos~\cite{sirignano2018dgm} and adapted to the MFG and PDE systems in~\cite{alaradi2018solving,carmona2021convergence,carmona2021deep}. The main motivation underlying this approach is to learn the PDE solutions using parameterized functions. This avoids computing the functions on a mesh, which is not feasible in high dimensions. In the DGM, we replace the function(s) solving the PDE(s) with a neural network(s), which are trained to minimize the PDE residual(s) as well as the boundary condition(s). 

To be specific, in our setting, we replace the functions $m$ and $u$ with neural networks, denoted by $m_{\theta}$ and $u_{\omega}$ and parameterized by $\theta$ and $\omega$ respectively. When the state $x$ is in high dimension, i.e., $d$ is large, we expect $m_{\theta}$ and $u_{\omega}$ to provide good approximations of $m$ and $u$ using much fewer parameters than the number of points in a grid. Furthermore, for the numerical implementation, we restrict our attention to a compact domain $\tilde\dom$. We denote $\tilde\domT = [0,T]\times\tilde\dom$. We expect the density to have a negligible mass outside a compact set so that by solving the PDE system on a large enough compact set, we obtain a good approximation of the solution, at least in the region where the density is significantly positive.  
We then define the loss function: %
\begin{equation*}
    \cL(\theta,\omega) = 
    \mathcal{L}^{({\mathrm{KFP}})}(m_{\theta},u_{\omega}) + \mathcal{L}^{({\mathrm{HJB}})}(m_{\theta},u_{\omega}),
\end{equation*}
where, for any $(m,u) \in \cC^{1,2}(\tilde\dom_T) \times \cC^{1,2}(\tilde\dom_T)$, the two losses are as:
\begin{align}
    \mathcal{L}^{({\mathrm{KFP}})}(m,u) 
    &= C^{({\mathrm{KFP}})}\left\| \displaystyle \frac{\partial m} {\partial t}  - \nu \Delta  m- \diver\Bigl( m\partial_p H(m,\nabla u)\Bigr) \right \| _{L^2(\tilde\dom_T)}^2 \nonumber
    \\ 
    &\qquad + C_{0}^{({\mathrm{KFP}})}\left\|m(0,\cdot) - m_{0} \right\|_{L^2(\tilde\dom)}^2
    + C_{T}^{({\mathrm{KFP}})}\left\|m(T,\cdot) - m_{T} \right\|_{L^2(\tilde\dom)}^2,
    \label{eq:DGM-loss-KFP}
\end{align}
and
\begin{align*}
    \mathcal{L}^{({\mathrm{HJB}})}(m,u) &=C^{({\mathrm{HJB}})}\Big\| \displaystyle
	-\frac{\partial u} {\partial t} - \nu \Delta u + H(m, \nabla u) 
    + \int_{\tilde\dom} \frac{\partial H} {\partial m}(\zeta, m(t,\cdot), \nabla u(t, \zeta))(\cdot) m(t,\zeta) \mathrm{d}\zeta \Big \|_{L^2(\tilde\dom_T)}^2.
\end{align*}
Here, $C^{({\mathrm{KFP}})}, C_{0}^{({\mathrm{KFP}})}, C_{T}^{({\mathrm{KFP}})}, C^{({\mathrm{HJB}})}$ are positive weights that give more or less importance to each component. If the space domain is bounded, we must include more penalty terms. Note that any smooth enough solution $(m,u)$ to the PDE system~\eqref{eq:MFOT-pde-system} makes $\mathcal{L}^{({\mathrm{KFP}})}$ and $\mathcal{L}^{({\mathrm{HJB}})}$ vanish. The goal is to find two neural networks which approximately minimize these losses.

Since it is not possible to compute exactly the above residuals, we approximate the $L^2$ norms using Monte Carlo samples. For example, we rewrite:
\begin{align*}
    &\left\| \displaystyle \frac{\partial m} {\partial t}  - \nu \Delta  m- \diver\Bigl( m\partial_p H(m,\grad u)\Bigr) \right \| _{L^2(\tilde \dom_T)}^2
    \\
    &= C(\tilde \dom_{T}) \cdot \mathbb{E}_{\tau,\xi} \left[ \left|\displaystyle \frac{\partial m} {\partial t}(\tau,\xi)  - \nu \Delta  m(\tau,\xi) - \diver\Bigl( m(\tau,\xi)\partial_p H(m(\tau),\nabla u(\tau,\xi))\Bigr) \right|^2 \right],
\end{align*}
where $(\tau,\xi)$ follows a uniform distribution over $\tilde \dom_T$, and $C(\tilde \dom_{T})$ is a normalizing constant that depends on the domain. Likewise, for the other norms, it would also be possible to use different norms and different distributions to sample $(\tau,\xi)$. But for the sake of simplicity, we will stick to this setting for the present work. We obtain the following probabilistic formulation of the loss function $\cL$:
\begin{equation*}
    \cL(\theta,\omega) = 
    \mathbb{E}_{S} \left[ \cL(\theta,\omega; S)\right], \qquad \cL(\theta,\omega; S) = \mathcal{L}^{({\mathrm{KFP}})}(m_{\theta},u_{\omega}; S) + \mathcal{L}^{({\mathrm{HJB}})}(m_{\theta},u_{\omega}; S), 
\end{equation*}
where $S = (\tau,\xi,\xi_0,\xi_T) \in [0,T] \times \tilde \dom \times \tilde \dom \times \tilde \dom$ denotes one sample, and for any $(m,u) \in \cC^{1,2}(\dom_T) \times \cC^{1,2}(\dom_T)$, the two losses at $S$ are as:
\begin{align*}
    \mathcal{L}^{({\mathrm{KFP}})}(m,u;S) 
    &= C^{({\mathrm{KFP}})} \left|\displaystyle \frac{\partial m} {\partial t}(\tau,\xi)  - \nu \Delta  m(\tau,\xi) - \diver\Bigl( m(\tau,\xi)\partial_p H(m(\tau),\nabla u(\tau,\xi))\Bigr) \right|^2
    \\ 
    &\qquad + C_{0}^{({\mathrm{KFP}})}|m(0,\xi_0) - m_{0}(\xi_0)|^2
    + C_{T}^{({\mathrm{KFP}})} |m(T,\xi_T) - m_{T}(\xi_T)|^2,
\end{align*}
and
\begin{align*}
    \mathcal{L}^{({\mathrm{HJB}})}(m,u;S) &=C^{({\mathrm{HJB}})}\Big| \displaystyle
	-\frac{\partial u} {\partial t}(\tau,\xi) - \nu \Delta u(\tau,\xi) + H(\xi, m(\tau), \nabla u(\tau,\xi)) 
    \\
    &\qquad\qquad\qquad + \int_{\dom} \frac{\partial H} {\partial m}(\zeta, m(t,\cdot), \nabla u(t, \zeta))(\xi) m(t,\zeta) \mathrm{d}y \Big|^2.
\end{align*}

Finally, to optimize over $(\theta,\omega)$, we use SGD (or one of its variants) on the loss $\cL$. In practice, we use a mini-batch of samples at each iteration, which amounts to approximate the expectation by an empirical average over several samples.

\subsection{Augmented Lagrangian Method with Deep Learning}
\label{subsec:deepadmm}

In this subsection, we present an approach based on a primal-dual formulation of the MFOT problem. We then introduce a deep learning adaptation of the alternating direction method of multipliers. We focus on the case when the interactions are local, and the drift is the control.

\subsubsection{Primal and dual problems}
Under suitable assumptions, the MFOT problem admits a variational formulation, which can be tackled using a direct optimization approach. 
As in the previous subsection, we assume that $\rho_0$ and $\rho_T$ have respectively density $m_{0}$ and $m_{T}$. %

We focus on a model with local interactions, meaning that an agent at state $x$ interacts with the density of the population at $x$. To alleviate the presentation, we will use the same notations for the costs and the drift functions, but now their second input is a real number $m$ instead of an element $\mu \in \cP_2(\dom)$. So we have $f: \dom \times \RR \times \RR^k \to \RR$, $g: \dom \times \RR \to \RR$, and $b: \dom \times \RR \times \RR^k \to \RR^d$. 
We also modify accordingly the definition of the Hamiltonian $H$ in \eqref{eq:def-hamiltonian} and the Lagrangian $L$ in \eqref{eq:def-lagrangian} in subsection~\ref{subsec:dgm}. We further assume that $f(x, m, v)$ is convex in $v$ for every $(x,m)$, and $mf(x, m, v)$ is convex in $m$ for every $(x,v)$. For simplicity, we consider that $b(x,m,v) = v$, i.e., the drift is the control. We remark that the setting here is not restrictive and satisfied by a large class of problems. \\

\textbf{Primal problem. } 
The MFOT problem~\eqref{eq:def-J-MFOT} introduced in Section \ref{sec:defpb} is formally equivalent to the following PDE-constrained optimization problem: 
\begin{align}
\label{eq:ADMM-original}
    & \inf_{v:\dom_{T} \rightarrow \RR^{k}} \, \int_{\dom_{T}} f\bigl(x, m(t,x), v(t,x)\bigr) m(x,t) \rd x \, \rd t \nonumber \\
    \text{subject to} \quad & \frac{\partial m}{\partial t}(t,x) - \nu \Delta m(t,x) + \diver \bigl( m(t, x) v(t, x) \bigr) = 0 \qquad t \in (0,T], x \in \dom \nonumber \\
    & m(0,x) = m_{0}(x), \qquad m(T,x) = m_{T}(x)
\end{align}
The PDE constraint is the KFP equation corresponding to the stochastic dynamics in~\eqref{eq:dyn-MFOT}. Note that the formulation in terms of $(m, v)$, while intuitive, is not convex in general. For this reason, we consider an equivalent formulation in terms of $(m,z) = (m, mv)$. We define:
\begin{align}
    \tilde{f}(x, m, z) 
    = \begin{cases}m f\left(x, m, \frac{z}{m}\right) & \text { if } m>0 \\ 0 & \text { if }(m, z)=(0,0) \\ +\infty & \text { otherwise }\end{cases}   
\end{align}
Note that $(m, z) \mapsto \tilde{f}(x,m,z)$ is LSC on $\RR \times \RR^{k}$. Under suitable conditions, it can be proved that $(m,z) \mapsto \tilde{f}(x,m,z)$ is convex on $\RR \times \RR^{d}$. We also define the space $\mathbf{K}$, 
\begin{align}
\label{eq:def-k1}
    \mathbf{K} = \Bigl\{(m,z) \,\Big|\,  \frac{\partial m}{\partial t}(t,x) - \nu \Delta m(t,x) + \diver z(t,x) = 0, m(0,x) = m_{0}(x), m(T,x) = m_{T}(x), m \geq 0 \Bigr\}
\end{align}
With all these definitions, we are ready to present the primal problem:
\begin{align}
\label{eq:ADMM-primal}
    \inf_{(m,z) \in \mathbf{K}} \cB(m,z) = \inf_{(m,z) \in \mathbf{K}} \int_{\dom_{T}} \tilde{f}\bigl(x, m(t,x), z(t,x) \bigr) \rd x \rd t
\end{align}
Assuming that problem~\eqref{eq:ADMM-primal} has a unique optimal solution $(m^*,z^*)$ and that problem~\eqref{eq:ADMM-original} has a unique optimal control $v^*$, then the following connection holds:
$v^{*}(t,x) = z^{*}(t,x)/m^{*}(t,x)$ if $m^{*}(t,x) > 0$, $v^{*}(t,x) = 0$ if $m^{*}(t,x) = 0$.

\textbf{Dual problem. } 
We now introduce a dual optimization problem. We define the following functionals:
\begin{align}
\label{eq:def-functional}
    & \cA(u) = \inf_{m \geq 0} \int_{\dom_{T}} m(t,x) \Bigl(\frac{\partial u}{\partial t}(t,x) + \nu \Delta u(t,x) - H\bigl(x, m(t,x), \nabla u(t,x) \bigr) \Bigr) \rd x \, \rd t 
    \\
    & \cF(u) = \int_{\dom} \left( m_{T}(x) u(T, x) - m_{0}(x) u(0,x) \right) \rd x 
    \\
    & \cG(\frak{a},\frak{b}) = - \inf_{m \geq 0} \int_{\dom_{T}} m(t,x) \Bigl(\frak{a}(t,x) - H\bigl(x, m(t,x), \frak{b}(t,x) \bigr) \Bigr) \rd x \, \rd t.
\end{align}
Note that if we define the linear differential operator $\Lambda u = \bigl(\frac{\partial u}{\partial t} + \nu \Delta u, \nabla u\bigr)$, then $\cA(u) = \cG(\Lambda u)$. Consider the following problem:
\begin{align}
\label{eq:ADMM-dual}
    \inf_{u} \cF(u) + \cG(\Lambda u).
\end{align}
Based on Fenchel-Rockafellar duality theorem (see Section 31, Theorem 31.1 in \cite{rockafellar1970convex}), we expect problems~\eqref{eq:ADMM-primal} and~\eqref{eq:ADMM-dual} to be in duality, meaning:
\begin{align}
\label{eq:primal-dual}
\inf_{(m,z) \in \mathbf{K}} \cB(m,z) = \sup_{u} \cA(u) = - \inf_{u} \cF(u) + \cG(\Lambda u) 
\end{align}
 Note that this primal-dual relationship also plays an important role in demonstrating the uniqueness and existence of solutions to MFG and MFC PDE systems, see e.g.~\cite{lasry2007mean,cardaliaguetgraber2015mean,achdou2016meancongestion1}. Here, we expect a similar result to hold for MFOT under suitable conditions. The rigorous definition of the two problems and the analysis of this duality relationship is left for future work. For now, we proceed formally.

We can at least formally establish a connection between the primal problem, the dual problem, and the optimal control in the following way. Let $u^{*}$ be the optimal solution to the dual \eqref{eq:ADMM-dual} and let $(m^{*}, z^{*})$ be the optimal solution to the primal problem~\eqref{eq:ADMM-primal}. Then the optimal control for the original problem~\eqref{eq:ADMM-original} is given by:
\begin{equation*}
    v^*(t,x) = \partial_p H\bigl(x, m^*(t,x),\nabla u^*(t,x)\bigr).
\end{equation*}

We notice that $(u^{*}, m^{*})$ forms a solution to the MFOT PDE system \eqref{eq:MFOT-pde-system}. This fact suggests that we can work on the dual problem \eqref{eq:ADMM-dual} directly to solve the MFOT problem. Under suitable assumptions, it can be shown that the dual problem \eqref{eq:ADMM-dual} is a strongly convex, unconstrained optimization problem over function space, which motivates the use of classic algorithms in convex optimization. However, the presence of the infinite dimensional linear operator $\Lambda$ makes the problem hard to solve efficiently in general. Fortunately, the structure of the objective as a sum of two convex functionals makes the problem amenable to algorithms based on splitting schemes, such as the Alternating Direction Method of Multipliers (ADMM) \cite{boyd2011distributed}. 

\begin{algorithm}
\DontPrintSemicolon
\KwData{Initial Guess $\bigl(u^{(0)}, q^{(0)}, \lambda^{(0)}\bigr)$; number of iterations $N$; hyperparameter $r > 0$}
\KwResult{Function $\bigl(u^{(N)}, q^{(N)}, \lambda^{(N)}\bigr)$ that are close to the saddle point of $\cL_r$ defined in~\eqref{eq:ADMM-augmented}}
\Begin{
\For{$k=1,\cdots,N, $}{
$u^{(k)} = \underset{u: \dom_{T} \rightarrow \RR}{\mathop{\arg\min}} \; \mathcal{F}(u) - \langle \lambda^{(k-1)}, \Lambda u \rangle + \frac{r}{2} \bigl\| \Lambda u - q^{(k-1)} \bigr\|^{2}$

$q^{(k)} = \underset{q: \dom_{T} \rightarrow \RR^{k+1}}{\mathop{\arg\min}} \; \mathcal{G}(q) + \langle \lambda^{(k-1)}, q \rangle + \frac{r}{2} \bigl\| \Lambda u^{(k)} - q\bigr\|^{2}$

$\lambda^{(k)} = \lambda^{(k-1)} - r\bigl(\Lambda u^{(k)} - q^{(k)}\bigr)$
}
}
\caption{Vanilla ADMM for MFOT} \label{algo:vanilla-admm}
\end{algorithm}

\subsubsection{Description of the algorithm}
Introducing a new variable $q$ that will play the role of $\Lambda u$, we can rewrite problem \eqref{eq:ADMM-dual} as the following constrained optimization program: 
\begin{align}
\label{eq:ADMM-false-constrained}
    \inf_{u, q:\, q = \Lambda u} \cF(u) + \cG(q).
\end{align}
The goal is now to find a saddle point of the associated Lagrangian. In fact, for numerical purposes, we will consider an augmented Lagrangian, defined as follows.

Let $r > 0$ be a constant and introduce $\lambda: \dom_{T} \rightarrow \mathbb{R}^{k+1}$, the Lagrangian multiplier associated with the constraint $q = \Lambda u$. Let $\inp{\cdot}{\cdot}$ denote the inner product on $L^{2}(\dom_{T})$. We introduce the augmented Lagrangian:
\begin{align}
\label{eq:ADMM-augmented}
    \mathcal{L}_r(u, q, \lambda) = \mathcal{F}(u) + \mathcal{G}(q) - \inp{\lambda}{\Lambda u - q} + \frac{r}{2} \| \Lambda u - q \|^{2}
\end{align}
Now, the original MFOT problem is reduced to finding a saddle point of $\cL_r$. Here, we state the original ADMM method in Algorithm \ref{algo:vanilla-admm} that finds the saddle point via an alternating optimization procedure. 

This general procedure can be implemented, for example, when the functions $(u, q, \lambda)$ are approximated by their values on a finite-difference grid. Such a procedure has been used for MFG and MFC problems, using finite elements~\cite{benamou2015augmented,andreev2017preconditioning} or finite differences~\cite{achdou2016meancongestion2}. Furthermore, \cite{benamou2015augmented} proved the convergence of this method under suitable conditions. However, as already mentioned, approximating functions by their values on a mesh is not feasible in high dimensions. We thus propose a different implementation of the ADMM based on neural network approximations.

In Algorithm \ref{algo:vanilla-admm}, the objectives in the steps are given by functionals to be minimized over functional spaces, which is not tractable in general. 
We restrict our attention to spaces of parameterized functions that can be expressed as neural networks, denoted by $\bigl(u_\theta, q_\omega, \lambda_\psi \bigr)$  with parameter $\theta, \omega, \psi$ respectively. We then follow the strategy introduced with the DGM \cite{sirignano2018dgm} and already used in Section~\ref{subsec:dgm} to create computable loss functions that are stochastic approximations of the functionals. 

Recall that the truncated space domain $\tilde \dom$ and the associated time-space domain $\tilde\dom$. Let $X \sim \mathcal{U}(\tilde\dom_{T})$, $Y \sim \mathcal{U}(\tilde\dom)$ be two random variables with uniform distribution in the time-space domain and the space domain respectively. Let $\rho_{X}$, $\rho_{Y}$ be the value of the uniform density on $\tilde \dom_{T}$ and $\tilde \dom$ respectively. %
Here, we overload the notation $\inp{\cdot}{\cdot}$ and $\| \cdot \|$ to represent the Euclidean inner product and norm on both $L_{2}(\tilde\dom_{T})$ and $\RR^{d+1}$:
$$
    \cL^{(u)}(\theta; \omega, \psi) = \cL_1(u_\theta,q_\omega,\lambda_\psi),
    \quad 
    \cL^{(q)}(\omega; \theta, \psi) = \cL_2(u_\theta,q_\omega,\lambda_\psi),
    \quad 
    \cL^{(\lambda)}(\psi; \theta, \omega, \psi_{old}) = \cL_3(u_\theta,q_\omega,\lambda_{\psi_{old}}, \lambda_\psi),
$$
where
\begin{align}
    \cL_1(u,q,\lambda) 
    & = \frac{1}{\rho_{Y}}\mathbb{E}_{Y} \left[ u(T, Y) m_{T}(Y) - u(0, Y) m_{0}(Y) \right] + \frac{1}{\rho_{X}} \mathbb{E}_{X} \left[\frac{r}{2}\|\Lambda u(X) - q(X) \|^{2} - \inp{\Lambda u(X)}{\lambda(X)} \right] %
    \label{eq:ADMM-loss-u}  
    \\
    \cL_2(u,q,\lambda) &= \frac{1}{\rho_{X}} \EE_{X} \left[ \cG(q(X)) + \inp{\lambda(X)}{q(X)} + \frac{r}{2} \|\Lambda u(X) - q(X) \|^{2} \right] \nonumber 
    \\
    \cL_3(u,q,\lambda_{old}, \lambda) &= \frac{1}{\rho_{X}} \EE_{X} \Bigl[ \|\lambda_{old}(X) - r \left( \Lambda u(X) - q(X)\right) - \lambda(X)\|^{2} \Bigr]. \label{eq:ADMM-loss-ql}
\end{align}
Here, the subscript $old$ is used to refer to the previous iteration: the loss for $\lambda$ involves the previous estimate $\lambda_{old}$. When using a neural network, it amounts to using the previous neural network parameters $\psi_{old}$. The loss function aims at mimicking the effect of the direct update in the third step of standard ADMM (Algorithm \ref{algo:vanilla-admm}) when $\lambda$ is approximated by a neural network. 

The algorithm DeepADMM is presented in Algorithm \ref{algo:DeepADMM}.

\begin{algorithm}
\DontPrintSemicolon
\KwData{Initial parameter $\theta^{(0)}, \omega^{(0)}, \psi^{(0)}$; number of ADMM iterations $K$; SGD parameters %
}
\KwResult{Final parameter $\theta^{(K)}, \omega^{(K)}, \psi^{(K)}$}
\Begin{
\For{$k=1,\cdots, K$}{
Compute $\theta^{(k)}$ using SGD to (approximately) minimize the loss $\cL^{(u)}(\cdot;\omega^{(k-1)},\psi^{(k-1)})$

Compute $\omega^{(k)}$ using SGD to (approximately) minimize the loss $\cL^{(q)}(\cdot;\theta^{(k)},\psi^{(k-1)})$

Compute $\theta^{(k)}$ using SGD to (approximately) minimize the loss $\cL^{(\lambda)}(\cdot;\theta^{(k)}, \omega^{(k)}, \psi^{(k-1)})$

}
}
\caption{DeepADMM for MFOT} \label{algo:DeepADMM}
\end{algorithm}

We have several remarks regarding DeepADMM and the augmented Lagrangian formulation in order for readers to better understand this approach. First, compared with Algorithm~\ref{algo:vanilla-admm}, the updates in Algorithm~\ref{algo:DeepADMM} for functions $u$ and $q$ are quite straightforward to understand. Instead of searching optimizer over function space, we reduce the problem to a finite dimension through stochastic approximation of the objective and search in the parameter space instead. The computed stochastic gradient can be considered as an unbiased estimator of the population gradient with respect to the functional, and the variance of this stochastic gradient decreases as the batch size increases. \\

\forlongversion{In Appendix~\ref{sec:ADMM-comp-G}, we discuss the computation of $\cG$ for several typical models.}

\section{Numerical experiments}
\label{sec:num-exp}

In this section, we present numerical experiments obtained with the three methods discussed in the previous section. For brevity, we refer to the three methods respectively introduced in sections~\ref{sec:method1}, \ref{subsec:dgm} and~\ref{subsec:deepadmm} as Method 1, Method 2 and Method 3 (and M1, M2, and M3 for short in the plots).

We first consider two test cases for which we have explicit solutions (up to solving ODE systems) and can thus be used to benchmark our algorithms in any dimension. We then consider two test cases that can be viewed as modifications of standard OT with crowd aversion or congestion effects.

\subsection{Case 1: Linear Quadratic Problem}
\label{sec:lq-test-case}

The first class of models that we consider has a linear-quadratic structure, which falls in the setting discussed in Example~\ref{ex:lq-example-description}. 

\subsubsection{Description of the problem}
In this model, we take: 
\begin{align*}
    b(x, \mu, a) = Ax + B a, \qquad 
    f(x, \mu, a) = a^{\top} R a, \qquad 
    \rho_0 = \cN(\bar x_0; \Sigma_0), \qquad 
    \rho_T = \cN(\bar x_T; \Sigma_T),
\end{align*}
where $A, B, R, \Sigma_0, \Sigma_T$ are (constant) matrices of suitable sizes. The vectors $\bar x_0$ and $\bar x_T$ correspond to the initial and terminal means.
We will consider two settings. In order to have a benchmark solution, we will take $\sigma = B$. This enables us to use the solution provided by~\cite[Section 7.1]{chen2018steering}, which boils down to solving a system of ODEs. \forlongversion{For the sake of completeness, we provide the details in Appendix~\ref{app:LQ-ODE}.}

\subsubsection{Evaluation Metrics}
\label{sec:evaluation-metrics}

In this model, since we have access to the optimal solution, we can evaluate the learnt solutions given by the three methods we proposed with respect to the ground-truth solution. We denote by $v^*$ the optimal control and $\hat v$ a learnt control. As explained below in detail, we use the following metrics: the total cost (namely $J^{MFOT}(\hat v)$, with $J^{MFOT}$ introduced in~\eqref{eq:def-J-MFOT}), the relative error between the achieved cost, and the optimal cost (namely $J^{MFOT}(\hat v)$ and $J^{MFOT}(v^{*})$), the deviation from the terminal distribution (i.e., the Wasserstein distance between the achieved terminal distribution and the target terminal distribution, $\rho_T$), and the weighted $L^{2}$ error between the learnt control $\hat v$ and the optimal control $v^{*}$, weighted by the population distribution.

\vskip 6pt

\noindent \textbf{Computation of the control.} The control is parameterized in different ways across different methods. For Method 1, $\hat v(t,x) = v_{\theta}(t,x)$. For Method 2 and Method 3, $\hat v(t,x) = -\frac{1}{2} B R^{-1} \nabla \hat u_\theta(t,x)$, where $\hat u_\theta$ is the neural network that approximates the dual variable, solution to the HJB equation.

\vskip 6pt

\noindent \textbf{Total cost.} Recall the definition of the objective $J^{MFOT}$ defined in \eqref{eq:def-J-MFC}. Let $v$ be a control. In the present Linear-Quadratic case, we have that,
\begin{align*}
    J^{MFOT}(v) = \int_{0}^{T} \int_{\cQ} m^{v}(t,x) f(x, m^{v}, v) \mathrm{d} x \, \mathrm{d} t = \int_{0}^{T} \int_{\cQ} m^{v}(t,x) v(t,x)^\top R v(t,x) \mathrm{d} x \, \mathrm{d} t,
\end{align*}
where $m^{v}$ is the density of mean field distribution driven by $v$, which satisfies the KFP PDE~\eqref{eq:MFOT-pde-system}.  In order to evaluate $J^{MFOT}(v)$, we use Monte Carlo simulations. We discretize the time variable $t$. Let $N_{T}$ be a number of time steps of length $\Delta t = T/N_{T}$. We consider a equi-distanced time discretization with time-steps $\{t_{0} = 0, t_{1} = \Delta t, \dots, t_{N_{T}} = N_{T} \Delta t\}$. Again, we simulate solutions to the underlying SDE using an Euler-Maruyama scheme similar to the one used in \eqref{eq:dyn-X-general-MFC-N-Deltat}. We simulate a family of $N$ sequences $((X_{t_{n}}^{i, v})_{n = 0, \dots, N_{T}})_{i = 1,\cdots,N})$ using the following update
\begin{equation}\label{eq:MC-simulation-lq}
\begin{cases}
    X^{i}_0 \sim \rho_0 \quad \hbox{ i.i.d.}
    \\
    X^{i,v}_{t_{n+1}}=X^{i,v}_{t_n} + (A X^{i, v}_{t_{n}} + B v(t_{n}, X^{i, v}_{t_{n}})) \Delta t + \sigma \sqrt{\Delta t} \Delta W^i_n,
\end{cases}
\end{equation}
where $\{\Delta W^{i}_{n}\}$ are i.i.d standard Gaussian random variables in $\mathbb{R}^{d}$. With these sampled sequences, we compute the objective as 
\begin{align*}
J^{MFOT}(v) = \dfrac{1}{N}\sum_{i=1}^N\sum_{n=0}^{N_T-1} v(t_n, X_{t_{n}}^{i, v})^{\top} R v(t_n, X_{t_{n}}^{i, v}) \Delta t.
\end{align*}

\vskip 6pt

\noindent \textbf{Relative Error.} The relative error between $J^{MFOT}(\hat v)$ and $J^{MFOT}(v^{*})$ is defined as 
$$
    \frac{|J^{MFOT}(\hat v) - J^{MFOT}(v^{*})|}{|J^{MFOT}(v^{*})|}.
$$ 
The major reason to consider relative error instead of absolute error is because the scale of the running cost varies greatly across different problems. Also, we want to stress that even though $v^{*}$ is the analytical optimal solution, it may happen that $J^{MFOT}(\hat v) < J^{MFOT}(v^{*})$ if $\hat v$ does not satisfy exactly the constraint (in contrast with $v^{*}$). %

\vskip 6pt

\noindent \textbf{Expected $L^{2}$ error for control.} The expected $L^{2}$ error between the learnt control $\hat v$ and the ground-truth control $v^{*}$ is defined as,
\begin{align*}
    d_{L^{2}}(\hat v, v^{*}) 
    = \int_{0}^{T} \int_{\dom} m^{*}(t,x) \|\hat v(t,x) - v^{*}(t,x) \|^{2} \mathrm{d}x \, \mathrm{d}t
\end{align*}
where $m^{*}$ is the density of the optimal mean-field associated with $v^{*}$. For the LQ problem, the optimal mean field $m^{*}$ is Gaussian for any $t \in [0, T]$, with mean $\mu_{t}$ and variance $\Sigma_{t}$ given by analytical formulas\forlongversion{ in Appendix~\ref{app:LQ-ODE}}. We can thus evaluate the $L^{2}$ error again with Monte Carlo samples for each time step.  As above, we discretize the time variable $t$ with $N_{T}+1$ points and for  $t_{n} = n\Delta t, n = 0, \dots, N_{T}$, we generate i.i.d. samples $(X^{i, v}_{t_{n}})_{i = 1, \dots, N} \sim \mathcal{N}(\mu_{t_{n}}, \Sigma_{t_{n}})$. We then estimate the $L^{2}$ error as: 
\begin{align*}
d_{L^{2}}(\hat v, v^{*}) \approx \dfrac{1}{N}\sum_{i=1}^N\sum_{n=0}^{N_T-1} \|\hat v(t_{n}, X_{t_{n}}^{i, v}) -  v^{*}(t_{n}, X_{t_{n}}^{i, v})\| ^{2}\Delta t .
\end{align*}

\vskip 6pt

\noindent \textbf{Deviation of distribution.} The deviation of the mean field $\hat \rho_{T}$ from the terminal target distribution $\rho_{T}$ is quantified by two different metrics: Wasserstein-2 distance $\cW_{2}(\hat \rho_{T}, \rho_{T})$ and $L^{2}$ distance $d_{L^{2}}(\hat m_{T}, m_{T})$. Here, $\hat \rho_{T}$ is the measure of the mean field distribution driven by the learnt control $\hat v$ at time $T$. $\hat m_{T}$ and $m_{T}$ are the density of $\hat \rho_{T}$ and $\rho_{T}$ respectively. 

\begin{itemize}
\item \textbf{Wasserstein-2 distance.} We adopt a similar method to compute the Wasserstein-2 distance as the one discussed in Section \ref{sec:method1}. We simulate $N$ particles following the dynamics~\eqref{eq:MC-simulation-lq}, and obtain a collection of $N$ samples $(X_{T}^{i, \hat v})_{i = 1, \dots, N}$, which forms an empirical distribution approximating $\hat \rho_{T}$. We also generate $N$ samples directly from the target distribution $\rho_{T}$, denoted by $(Y^{i}_{T})_{i = 1, \dots, N}$. Then, we define the distance matrix $\cM$ by $\cM_{ij} = |X_{T}^{i , \hat v} - Y_{T}^{j}|^{2}$, and we recall that the set $U_{N}$ is defined by~\eqref{eq:set-U-Wasserstein}.  
We approximate the Wasserstein-2 distance between the two empirical distributions formed by $X_{T}^{i, 
\hat v}$ and $Y_{T}^{i}$ through the following linear program:
\begin{align*}
    \mathcal{W}_2\left(\hat \rho_{T}, \rho_{T} \right) \approx \Bigl(\min\limits_{\cA \in U_N} \inp{\cA}{M}\Bigr)^{1/2}.
\end{align*}
\item \textbf{$L^{2}$ distance.} We will also use as a metric the $L^{2}$ distance between $\hat m_{T}$ and $m_{T}$ on the truncated domain $\tilde\dom$.
To evaluate the integral, in the absence of analytical formula for $m_T$, we again use a Monte Carlo approach. We uniformly sample $N$ points in $\tilde\dom$ denoted by $(X_{j})_{j = 1, \dots, N}$. Let $C(\tilde\dom)$ denotes the inverse of the value of the uniform density. Then we can approximate the $L^{2}$ distance by:
\begin{align*}
     \frac{C(\tilde\dom)}{N} \sum_{i = 1}^{N} \|\hat m_{T}(X_{i}) - m_{T}(X_{i}) \|^{2}. 
\end{align*}

\end{itemize}

\subsubsection{Numerical results}

In the numerical tests, we take the values given in Table~\ref{table:LQ-params} for the parameters of the model, with the time horizon $T=1.0$. 
\begin{table}[t]
\begin{center}
\begin{tabular}{||c | c | c | c | c | c | c | c | c | c ||} 
 \hline
 Test & $d$ & $A$ & $B$ & $\sigma$ & $R$ & $\bar x_0$ & $\Sigma_0$ & $\bar x_T$ & $\Sigma_T$ \\ 
 \hline
 \tabincell{c}{LQ Test 1} & $1$ & $1$ & $1$ & $1$ & $\frac{1}{2}$ & $0.0$ & $1$ & $2.0$ & $0.5$
 \\
 \hline
 \tabincell{c}{LQ Test 2} & $2$ & $I_d$ & $I_d$ & $I_d$ & $\frac{1}{2}I_d$ & $[0.0, 0.0]$ & $I_d$ & $[2.0, 2.0]$ & $\frac{1}{2}I_d$
 \\
 \hline
\end{tabular}
\caption{Parameters for the two linear-quadratic test cases}
\label{table:LQ-params}
\end{center}
\end{table}

\begin{figure}[t]
    \centering
    \includegraphics[width=0.98\textwidth]{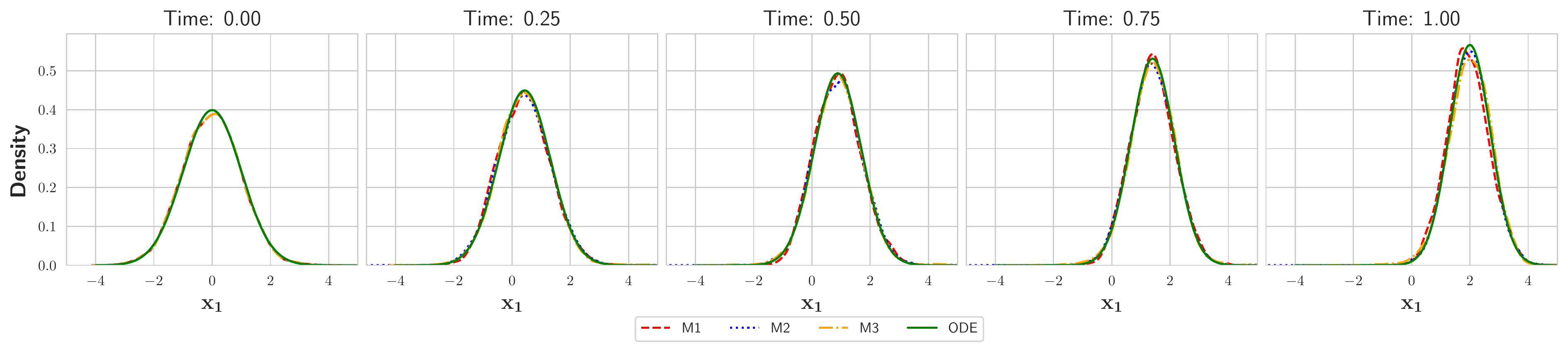}
    \caption{Evolution of the density in the LQ Test case 1. Each plot corresponds to one time step and displays the densities as functions of the space variable, $x$. The densities are: The density obtained by applying the control learnt by each of the three deep learning methods as well as the ground-truth density given by the ODE method.}
    \label{fig:lq_test1_evol}
\end{figure}

\begin{figure}[t]
    \centering
    \includegraphics[width=0.98\textwidth]{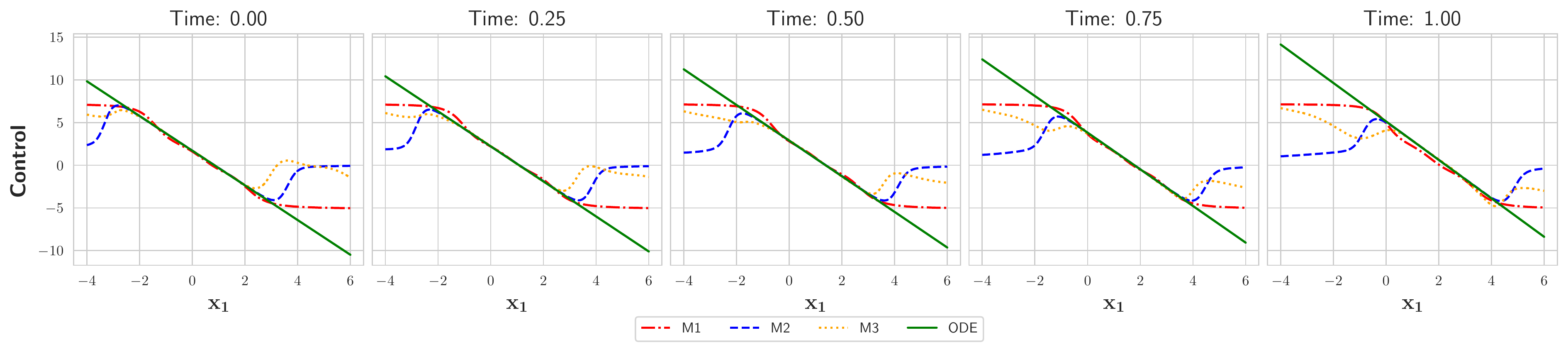}
    \caption{Evolution of the control in the LQ Test case 1. Each plot corresponds to one time step and displays the controls as functions of the space variable, $x$. The controls are: The control learnt by each of the three deep learning methods as well as the ground-truth control given by the ODE method.}
    \label{fig:lq_test1_control}
\end{figure}

For LQ test 1, in dimension 1, Figure~\ref{fig:lq_test1_evol} displays the evolution of the density. For methods 1, 2, and 3, we obtain the control learnt using the neural networks and simulate $N$ trajectories by the Monte Carlo method following the dynamics~\eqref{eq:MC-simulation-lq}, with $v$ replaced by the learnt control. We then estimate the mean field distribution using kernel density estimation (KDE). We see that the distributions obtained with the three methods match well the ground-truth one obtained with ODEs. The distributions move towards the right and concentrate around the final mean. Figure~\ref{fig:lq_test1_control} shows the evolution of the control. We see that the three methods provide good approximations of the true optimal control, at least in the region where the density is high. In regions where the density is very low, the control is not well approximated, but this is not an issue as far as the optimal behavior of the population is concerned. The first part of Table~\ref{tb:lq_test1-result} shows the results obtained for the metrics introduced above. We see that each of the three methods achieves a smaller total cost than the true optimal control. This would not be possible for controls satisfying perfectly the terminal constraint, but it is possible here due to the fact that the methods satisfy only approximately the planning constraint. The optimal control is well approximated, as shown by the $L^2$ distance to the true optimal control. Furthermore, we see that Methods~2 and~3 have a higher Wasserstein-2 distance between the terminal distribution and the target distribution, but the $L^2$ distance is much lower. %

\begin{figure}[t]
    \centering
    \includegraphics[width=0.8\textwidth]{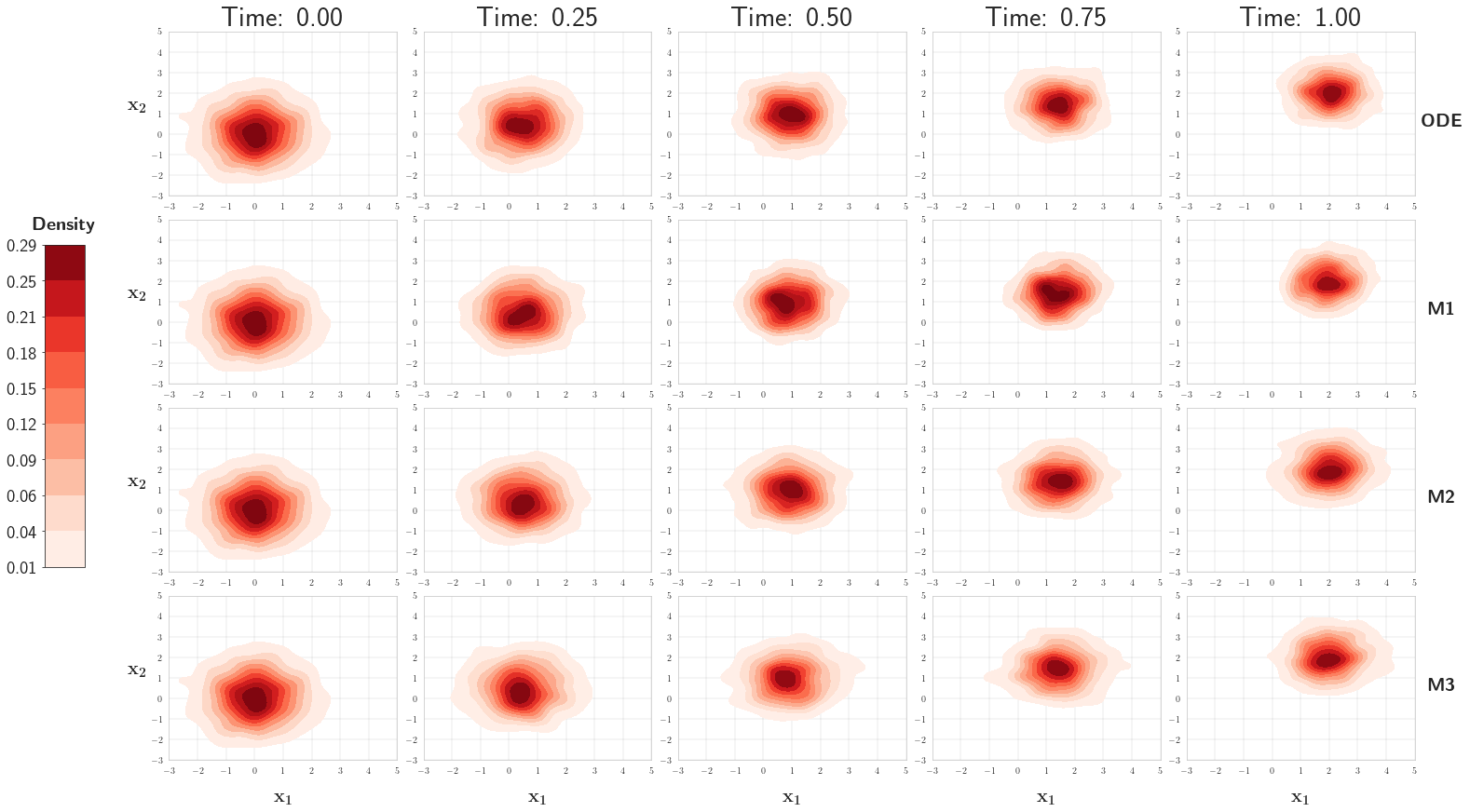}
    \caption{Evolution of the density in the Linear Quadratic Test case 2. Each column corresponds to one time step, and each row corresponds to one of the methods. Each plot displays the density as a function of the space variable, i.e., $(m(t,x)_1)_{x \in [-4,6]^2}$. The first row corresponds to the solution obtained by the ground-truth ODE method. The second, third and fourth rows correspond respectively to methods 1, 2 and 3.}
    \label{fig:lq_test2_evol}
\end{figure}

\begin{figure}[h]
    \centering
    \includegraphics[width=0.8\textwidth]{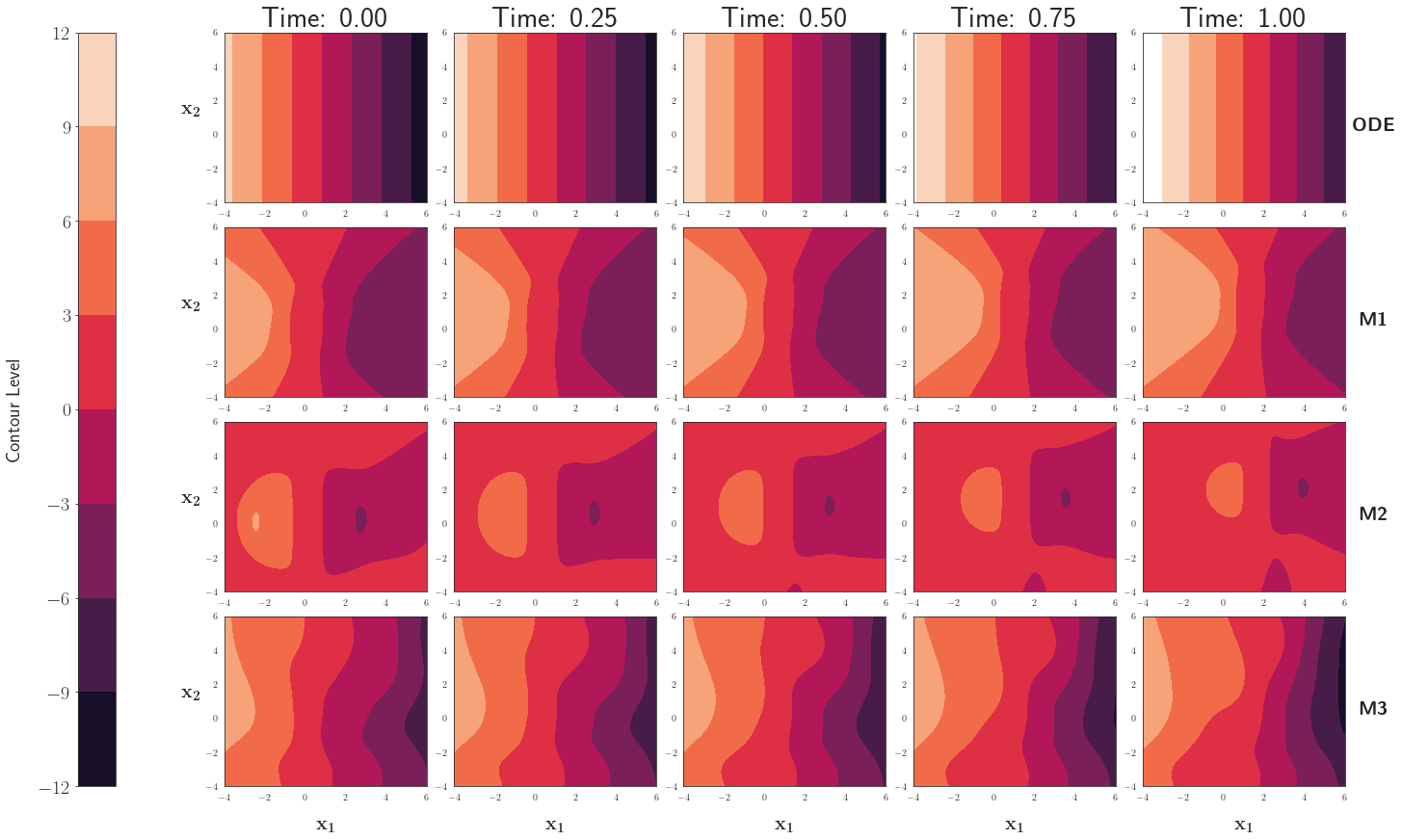}
    \caption{Evolution of the control in the Linear Quadratic Test case 2. Each column corresponds to one time step, and each row corresponds to one of the methods. Each plot displays the first dimension of the control as a function of the space variable, i.e., $(v(t,x)_1)_{x \in [-4,6]^2}$. The first row corresponds to the solution obtained by the ground-truth ODE method. The second, third and fourth rows correspond respectively to methods 1, 2 and 3.}
    \label{fig:lq_test2_control_0}
\end{figure}

As for LQ test 2, in dimension 2, Figure~\ref{fig:lq_test2_evol} displays the evolution of the density for each of the methods. The densities move from the bottom left corner to the top right corner. Furthermore, the terminal distribution is more concentrated because the terminal variance is smaller than the initial variance. Figure~\ref{fig:lq_test2_control_0} shows the evolution of the first dimension of the control (the second dimension is similar, so we omit it for brevity). The ground-truth control is linear in space for each time step. We see that the three methods manage to learn approximately linear controls, at least in the region where the density is significantly positive. 
Table~\ref{tb:lq_test2-result} shows the results obtained for the metrics introduced above. We see that here again, each of the three methods achieves a smaller total cost than the true optimal control due to the fact that the terminal constraint is not perfectly satisfied. The optimal control is well approximated, and the terminal distribution is matched with good accuracy.

\begin{table}[t]
\begin{center}
\begin{small}
\begin{tabular}{||c|c|c|c|c|c|c|c||}
\hline
Test case & Method & Total Cost & Relative Error & $d_{L^{2}}(\hat v, v^{*})$ & $\mathcal{W}_{2}(\hat \rho_{T}, \rho_{T})$ & $d_{L^{2}}(\hat m_{T}, m_{T})$\\
\hline
\multirow{4}{*}{\tabincell{c}{Linear Quadratic \\ LQ Test 1}} 
& ODE $(v^{*})$& $2.126$ & \,\, - \,\, &  \,\, - \,\,  & \,\, - \,\, & \,\, - \,\,  \\
& M1 & $2.099$ & $1.24\%$ & $0.021$ & $0.002$ & $0.006$ \\
& M2 & $2.096$ & $1.41\%$ & $0.003$ & $0.043$ & $0.00004$ \\
& M3 & $2.077$ & $2.29\%$ & $0.011$ & $0.031$ & $0.001$\\
\hline
\end{tabular}
\caption{Comparison of three different methods v.s. the analytical solution on the LQ test case 1. The evaluation metrics are described in Section~\ref{sec:evaluation-metrics}.}
\label{tb:lq_test1-result}
\end{small}
\end{center}
\end{table}

\begin{table}[t]
\begin{center}
\begin{small}
\begin{tabular}{||c|c|c|c|c|c|c|c||}
\hline
Test case & Method & Total Cost & Relative Error & $d_{L^{2}}(\hat v, v^{*})$ & $\mathcal{W}_{2}(\hat \rho_{T}, \rho_{T})$ & $d_{L^{2}}(\hat m_{T}, m_{T})$\\
\hline
\multirow{4}{*}{\tabincell{c}{Linear Quadratic \\ LQ Test 2}} 
& ODE $(v^{*})$ & $4.175$ & \,\, - \,\, &  \,\, - \,\,  & \,\, - \,\, & \,\, - \,\,  \\
& M1 & $4.117$ & $1.39\%$ & $0.074$ & $0.043$ & $0.007$\\
& M2 & $3.935$ & $5.63\%$ & $0.043$ & $0.403$ & $0.00005$\\
& M3 & $4.054$ & $2.89\%$ & $0.131$ & $0.561$ & $0.015$\\
\hline
\end{tabular}
\caption{Comparison of three different methods v.s. the analytical solution on the LQ test case 2. The evaluation metrics are described in Section~\ref{sec:evaluation-metrics}.}
\label{tb:lq_test2-result}
\end{small}
\end{center}
\end{table}

\subsection{Case 2: Transport with congestion effects}

The second class of models that we consider is inspired by crowd motion and falls in the setting discussed in Example~\ref{ex:congestion-example-description}.

\subsubsection{Description of the problem}

In this model, intuitively, the cost is higher when moving through a crowded region, i.e., where the density is high. Specifically, we take:
\begin{align*}
    b(x, \mu, a) = a, \qquad 
    f(x, \mu, a) = R |\ell(x,\mu)|^{\gamma} |a|^2, \qquad 
    \rho_0 = \cN(\bar x_0; \Sigma_0), \qquad 
    \rho_T = \cN(\bar x_T; \Sigma_T). 
\end{align*}
For the function $\ell$, we take two different models. We consider the following non-local dependence: 
$$
    \ell(x, \mu) = c + \rho_\epsilon \star \mu(x),
$$
where $c>0$ is a constant, $\rho_\epsilon$ is a Gaussian kernel and $\star$ denotes the convolution. We use Method 1 to solve the MFOT problem with this function $\ell$. Since it is based on Monte Carlo simulations of trajectories, it is straightforward to compute a convolution with the empirical distribution at a given time step. 

We also consider a variation with a local dependence.
$$
    \ell(x, \mu) = c + m(x)
$$
where $c>0$ is a constant and $m$ denotes the density of $\mu$. For this type of model, Methods 2 and 3 are better suited since, in these methods, we directly have access to the approximate density in the form of a neural network.

\subsubsection{Numerical results} 
\label{subsec:cong_num}

We focus on one test case called "Congestion" below in dimension $d=1$. In this model, $\gamma=1$. For the sake of comparison, we also consider the corresponding model with the same choice of parameters except that $\gamma=0$, i.e., there are no congestion effects in the running cost. The values that we take in the numerical tests are given in Table \ref{tab:cong_params} below. 
\begin{table}[h]
\begin{center}
\begin{tabular}{||c | c | c | c | c | c | c | c | c | c||} 
 \hline
 Test case & $d$ & $\gamma$ & $c$ & $\sigma$ & $R$ & $\bar x_0$ & $\Sigma_0$ & $\bar x_T$ & $\Sigma_T$ \\ %
 \hline
 Case 1, No congestion& $1$ & $0$ & $0.1$ & $0.1$ & $0.5$ & $0$ & $0.04$ & $2$ & $0.04$ 
 \\
 \hline
 Case 2, Congestion & $1$& $1$ & $0.1$  & $0.1$ & $0.5$ & $0$ & $0.04$ & $2$ & $0.04$ 
 \\
 \hline
 Case 3, Congestion& $5$&  $1$ & $1$  & $1$ & $0.5$ & $0$ & $0.1$ & $2$ & $0.1$ 
 \\
 \hline
\end{tabular}
\caption{Parameters for the test case with congestion and the benchmark model without congestion effects}
\label{tab:cong_params}
\end{center}
\end{table}

\begin{figure}[t]
    \centering
    \includegraphics[width=0.98\textwidth]{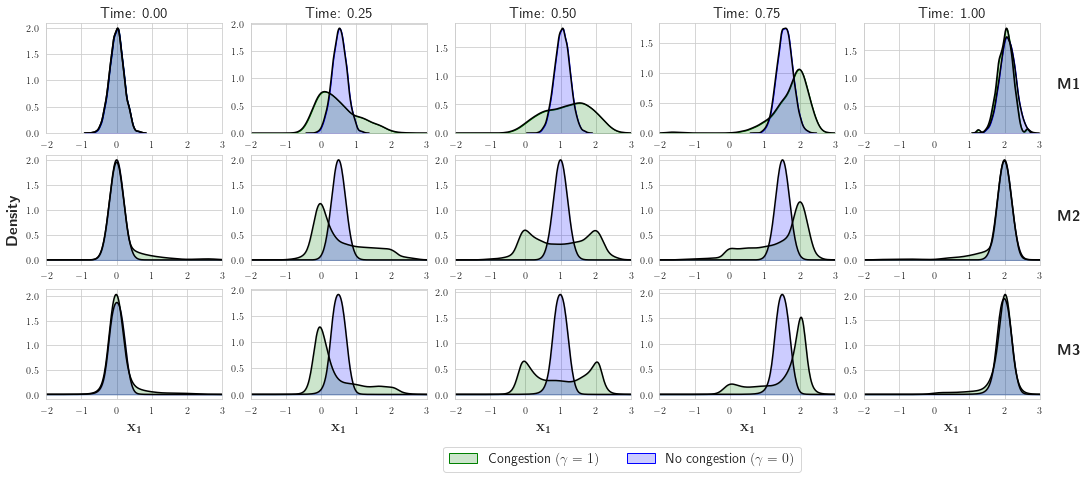}
    \caption{Visualization of the mean field density $\hat m(t,x)$ in Congestion test Case 1,2}
    \label{fig:cong_test1_evol}
\end{figure}

In Figure~\ref{fig:cong_test1_evol}, we present the evolution of the density under the control learnt by each of the three methods for congestion cases 1 and 2. Each row corresponds to one method. We see that, in the case where $\gamma=0$ (no congestion effect), the mass is transported directly towards the terminal distribution without much change in its shape. In contrast, in the case with $\gamma=1$, the mass spreads in space and one part starts moving towards the target mean $\bar x_T = 2$ whereas another part stays behind and catches up at later time steps. This is consistent with the idea that moving in congested regions is more expensive, so some agents would agree to wait until the density decreases before moving forward.

\begin{figure}[t]
    \centering
    \includegraphics[width=0.98\textwidth]{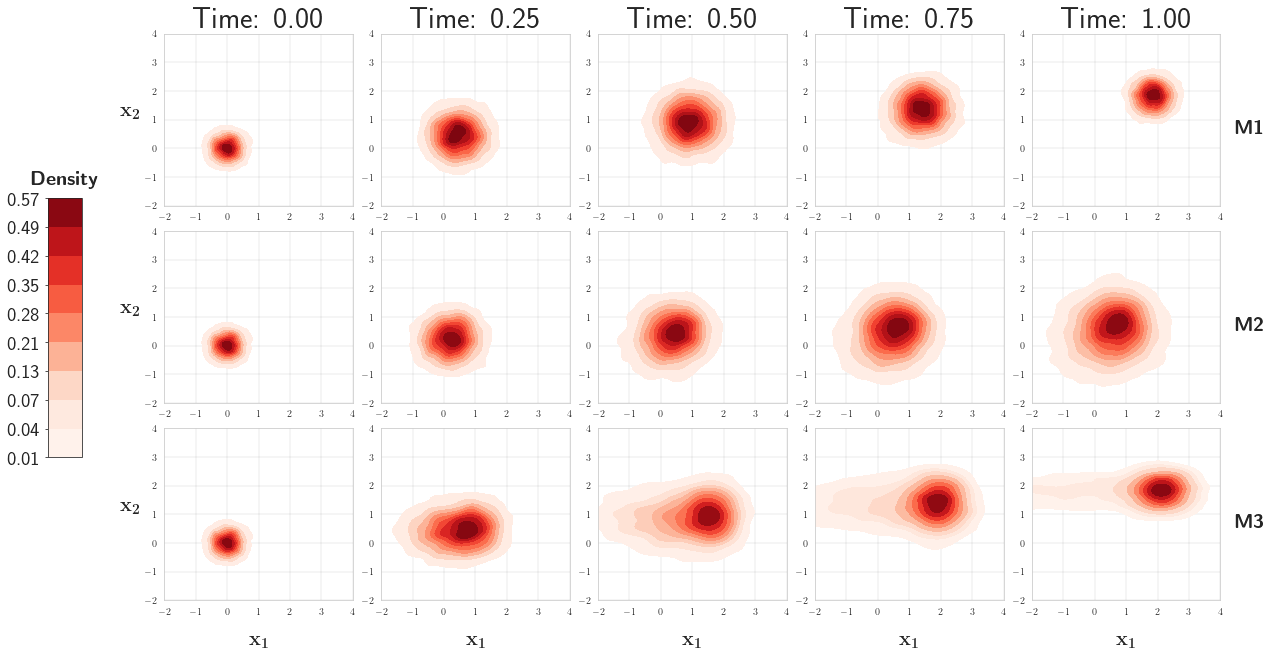}
    \caption{Visualization of the mean field density $\hat m(t,x)$ in Congestion test Case 3}
    \label{fig:cong_test3_evol}
\end{figure}

Finally, in Figure~\ref{fig:cong_test3_evol}, we present the evolution of the density under the control learnt by each of the three methods for congestion case 3, which is in dimension 5. Each row corresponds to one method. To visualize density evolution in dimension 5, we plot the marginal distribution of the mean field distribution on the first and second dimensions. We see that, similarly to the congestion case 2, the mass spreads in space and gradually moves towards the target mean. Compared with congestion case 2, the difference in the moving pattern and extent of spreading is due to the difference of parameters in Table~\ref{tab:cong_params}. With a larger value $c$, the behavior of the density would be closer to a direct transport to the terminal distribution without changes in the shape of the distribution.

\subsection{Remarks on the choice of hyperparameters}

Each method has several hyperparameters, including the architecture of the neural networks. We provide below some remarks about the choice of hyperparameters in our implementation. 

\textbf{Method 1. } In our implementation, we choose $G(r) = C_W r$ where $C_W$ is a hyperparameter that we adjust dynamically. We increase the constant $C_W$ when we expect a higher running cost (for instance, in a higher dimension) in order to give enough importance to the penalty. The coefficient $\alpha$ of regularization for the computation of the Wasserstein distance is also a hyperparameter that we adjust dynamically using the following heuristics. We start with a given value for $\alpha$ and, when the estimated Wasserstein distance is small enough, we reduce the value of $\alpha$. The idea is that, as long as the terminal distribution does not match well enough the target distribution, we need a high level of regularization in order to estimate efficiently the Wasserstein distance between them. As the two distributions get closer, we can decrease the degree of regularization in order to have a more accurate estimation of the Wasserstein distance. The way we adjust $C_W$ also depends on the dimension of the state variable. There is also a computational time aspect to take into account: as $\alpha$ becomes smaller, the computations take more time (see~\ref{subsec:method1-algo} for more details). For the neural network, we take a feedforward fully connected neural network with 6 layers of 60 neurons each. The other hyperparameters are the number of particles $N$ and the number of time steps $N_T$. We take $N=300$ and $N_T=20$.

\vskip 6pt

\textbf{Method 2. } In the second method, no time or space discretization is needed, and the density is directly approximated by a neural network, so we do not need to use a finite number of particles. However, we need to choose the values of the weights $C_{0}^{({\mathrm{KFP}})},$ $C_{T}^{({\mathrm{KFP}})},$ $C^{({\mathrm{KFP}})},$ and $C^{({\mathrm{HJB}})}$ in the loss function. We used $C_{0}^{({\mathrm{KFP}})} = 20,$ $C_{T}^{({\mathrm{KFP}})} = 50,$ $C^{({\mathrm{KFP}})} = 20,$ $C^{({\mathrm{HJB}})} = 1$. As for the neural network, we used the architecture proposed in the DGM article~\cite{sirignano2018dgm}, with 2 layers and a width equal to 40.  During the training, at each iteration of SGD, we use a minibatch of 500 points in time and space, and 500 points in space for the initial and terminal conditions. %

\vskip 6pt

\textbf{Method 3. } The main hyperparameter in this method is $r$, which is used in the definition of the augmented Lagrangian~\eqref{eq:ADMM-augmented}. For the experiments, we select $r = 0.1$. Even though, in theory, the convergence of ADMM is independent of the choice of $r$, in practice, we often find that a large $r$ value could potentially increase numerical instability and lead the algorithm to diverge. Similarly, a small $r$ value could slow down the convergence. As for the neural networks, we use the following architectures. For both $u_{\theta}$, $q_{\omega}$, and $\lambda_{\psi}$, in general, we use a fully connected neural network with residual connections, sigmoid activation function, and appropriate output dimension. We use $6$ layers and $100$ neurons per layer. For LQ test cases, we further consider an extra quadratic correction in addition to the neural networks: the output of $u_{\theta}$ is the sum of neural network output and a quadratic function with trainable weights. To effectively model the mean field density, a sigmoid activation function is applied to the first dimension of the output of the neural network $\lambda_{\psi}$, and then the result is multiplied by a constant $C$. In this way, the first dimension of $\lambda_{\psi}$ takes values in $(0,C)$. In the experiments, we take $C = 1$ for the LQ test cases and $C = 5$ for the congestion test cases. During training, at each iteration of SGD, we use a minibatch of 512 points in time and space, and 512 points in space for the initial and terminal conditions.

\section{Conclusion and future directions}

In this work, we have proposed three numerical methods based on deep learning for mean field optimal transport problems. The three methods can tackle a larger class of problems than deep learning methods proposed previously, which were mostly focusing on the Schr\"odinger bridge problem or MFGs with a specific structure. The first method replaces the terminal constraint with a penalty and then directly learns the optimal control using Monte Carlo trajectories. The second method solves a PDE system which is obtained as the optimality conditions for the MFOT problem. The third method relies on an augmented Lagrangian approach for the variational formulation of the problem. The numerical results show that the three methods match the analytical solution on an LQ problem, and that they are able to handle non-trivial mean field interactions modeling congestion effects.

From here, we can envision several research directions. First of all, the theoretical analysis of the MFOT problem remains to be tackled. For example, the existence and uniqueness of the solution to the PDE system have been proved only in relatively specific cases, see e.g.~\cite{achdou2012mean,cardaliaguet2013geodesics,porretta2013planning,graber2019planning}. It would be interesting to extend the analysis to more general forms of dynamics and cost functions. From the numerical point of view, it would be interesting to scale-up the methods proposed in this work to a higher dimension, and to explore other deep learning methods. The numerical analysis and the convergence proof of the proposed methods also remain to be investigated in future work.

\bibliographystyle{amsplain}
\bibliography{bibmfg}

\clearpage

\appendix

\section{Solution for LQ problem}
\label{app:LQ-ODE}
In this section, we provide an explicit solution for the LQ setting considered in Section~\ref{sec:lq-test-case}. We summarize the analytical solution derived in~\cite[Section 7.1]{chen2018steering} and we add the analytical formula for the variance, which is useful to fully describe the evolution of the population distribution. 
We consider the following setting, which is actually slightly more general than the one used in Section~\ref{sec:lq-test-case}. For any probability distribution $\mu$ admitting a first moment, we use the notation $\bar\mu = \int x \mu(dx)$, to be understood coordinate-wise. In this model, we take: 
\begin{align*}
    b(x, \mu, a) = Ax + \bar A \bar \mu + B a, \qquad \sigma = B, \qquad 
    f(x, \mu, a) = \frac{1}{2} a^{\top} a, \qquad 
    \rho_0 = \cN(\bar x_0; \Sigma_0), \qquad 
    \rho_T = \cN(\bar x_T; \Sigma_T).
\end{align*}
We denote by $v^*$ the optimal control and by $(\mu_t)_{t \in [0,T]}$ the optimal flow of distributions. Since the dynamics are linear and the initial and terminal distributions are Gaussians, it can be shown that the distribution remains Gaussian at every $t$. At each $t$, we denote by $\bar\mu_t$ its mean and by $\Sigma_t$ its covariance matrix.  
In the case without mean field interactions, the optimal control is linear in space, and the coefficient depends on the solution of a Ricatti equation in time. The mean field interaction being uniform in space and the fact that our problem is linear explains why we can set 
$$v^*(t,x)=- B^\top \Pi_t x + B^\top n_t$$
and solve for $(\Pi_t)_{t\in[0,T]}$ and $(n_t)_{t\in[0,T]}$, as time-dependent functions only. Note that $\bar \mu$ is now solution of the following ODE depending on $n_t$:
\begin{equation}
\label{eq:LQ2-MFOT-expli-mean-brieuc}
    \dot{\bar \mu}_t = (A + \bar A - BB^{\top} \Pi_{t}) \bar\mu_t + B B^\top n_t, \quad \bar \mu_0 = \bar x_0. %
\end{equation}
Furthermore, if we apply It\^o's formula to $V_t=\mathbb{E}[X_{t} X_{t}^{\top}]=\Sigma_t+\bar \mu_t \bar\mu_t^\top$, we also obtain a matrix ODE depending on $(n_t)_{t\in[0,T]}$ for $(V_t)_{t\in[0,T]}$:
\begin{align}
\label{eq:lq-variance-brieuc}
\begin{cases}
    \dot V_{t} = V_{t}(A^{\top} - \Pi_{t}^{\top}BB^{\top}) + (A - BB^{\top}\Pi_{t})V_{t} + \bar \mu_{t} n_{t}^{\top}BB^{\top} + BB^{\top}n_{t}\bar \mu_{t}^{\top} + \bar \mu_{t} \bar \mu_{t}^{\top} \bar A^{\top} + \bar A \bar \mu_{t} \bar \mu_{t}^{\top} + BB^{\top}
    \\
    V_{0} = \Sigma_{0} + \bar \mu_{0} \bar \mu_{0}^{\top}.
\end{cases}
\end{align}
We will also need the state transition matrices associated to $A$ and $A+\bar A$: we define $\Phi$ as the solution of 
\begin{equation*}
    \frac{\partial\Phi(t,s)}{\partial t} = A \Phi(t,s), \quad  \Phi(s,s) = I, \quad 0 \le s \le t \le T.
\end{equation*}
Since $A$ is constant, we have $\Phi(t,s)=e^{(t-s)A}$. Similarly, we define $\bar\Phi(t,s)=e^{(t-s)(A+\bar A)}$. We also set 
\begin{align} 
\begin{cases}
    M(t,s) = \displaystyle\int_s^t \Phi(t,\tau) B B^\top \Phi(t,\tau)^\top \mathrm{d}\tau, \quad 0 \le s \le t \le T,
    \\
    \bar M(t,s) = \displaystyle\int_s^t \bar\Phi(t,\tau) B B^\top \bar\Phi(t,\tau)^\top d\tau, \quad 0 \le s \le t \le T,
\end{cases}
\end{align}
and denote for brevity $\Phi_{T,0}=\Phi(T,0),$ $\bar \Phi_{T,0}=\bar \Phi(T,0),$ $M_{T,0}=M(T,0),$  and $\bar M_{T,0}=\bar M(T,0)$. 
Notice that if $BB^{\top}$ commutes with $\Phi$ and $\bar \Phi$, and if the matrices $A + A^{\top}$ and $A + \bar A + A^{\top} + \bar A^{\top}$ are non-singular, then we have the following closed form expressions for $M(t,s)$ and $\bar M(t,s)$, with $0 \le s \le t \le T$:
\begin{align*}
    & M(t,s) = BB^{\top} \int_{s}^{t} e^{A + A^{\top}}(t - \tau) \mathrm{d} \tau = BB^{\top} \bigl[e^{(A + A^{\top}) t} - e^{(A + A^{\top}) s}\bigr] (A + A^{\top})^{-1} 
    \\
    & \bar M(t,s) = BB^{\top} \int_{s}^{t} e^{A + A^{\top} + \bar A + \bar A^{\top}}(t - \tau) \mathrm{d} \tau = BB^{\top} \bigl[e^{(A + A^{\top} + \bar A + \bar A^{\top})t} - e^{(A + A^{\top} + \bar A + \bar A^{\top})s}\bigr] (A + A^{\top} + \bar A + \bar A^{\top})^{-1}.
\end{align*}

Now, $(\Pi_{t})_{t \in [0,T]}$ is the solution to the following Riccati ODE, which is independent from the other variables: %
\begin{align} 
\label{eq:lq-riccati}
\begin{cases}
    \dot\Pi_t = - A^\top \Pi_t - \Pi_t A + \Pi_t B B^\top \Pi_t
    \\
    \Pi_0 = \Sigma_0^{-1/2}\left[ \frac{I}{2} + \Sigma_0^{1/2}\Phi_{T,0}^\top M_{T,0}^{-1}\Phi_{T,0} \Sigma_0^{1/2} - \left( \frac{I}{4} + \Sigma_0^{1/2}\Phi_{T,0}^\top M_{T,0}^{-1} \Sigma_T M_{T,0}^{-1}\Phi_{T,0} \Sigma_0^{1/2} \right)^{1/2} \right]\Sigma_0^{-1/2}.
\end{cases}
\end{align}
Under some conditions, Riccati equations admit explicit solutions in dimension $1$, see e.g. page 110 in~\cite{carmona2018probabilistic}. More generally, we can solve~\eqref{eq:lq-riccati} using a forward time-marching method. 

We can then show that  %
\begin{equation}
\label{eq:nt-lq-with-mf}  
    n_t = \Pi_t \bar\Phi(t,T) \bar M(T,t) \bar M_{T,0}^{-1} \bar \Phi_{T,0} z_0 + \Pi_t \bar M(T,0) \bar\Phi(T,t)^\top \bar M_{T,0}^{-1} z_T + \bar \Phi(T,t)^\top \bar M_{T,0}^{-1}(z_T - \bar \Phi_{T,0} z_0),
\end{equation}
where $(z_t)_{t\in[0,T]}$ is defined by
\begin{equation}
\label{eq:mt-lq-with-mf}    
    z_t = \bar\Phi(T,t)^\top \bar M_{T,0}^{-1} (\bar x_T - \bar\Phi_{T,0} \bar x_0).
\end{equation}

Now, up to the computation of $\Pi$, we have an explicit formula for $n$ and we can obtain the mean $\bar\mu_t$ from \eqref{eq:LQ2-MFOT-expli-mean-brieuc} and then the covariance matrix $\Sigma_t$ from \eqref{eq:lq-variance-brieuc}. The optimal mean field distribution at time $t$ is the Gaussian distribution with mean $\bar \mu_{t}$ and variance $\Sigma_{t}$.

\section{Computation of $\cG$}
\label{sec:ADMM-comp-G}

In this section, we discuss some of the issues arising when computing $\cG$ in practice, as well as our solutions. As the reader may notice, the losses defined in \eqref{eq:ADMM-loss-u} and \eqref{eq:ADMM-loss-ql} depend on the exact form of functionals $\cF$ and $\cG$. As $\cF$ is already defined in an explicit, easy-to-compute form, we are left with the problem of figuring out a good approach to compute $\cG$. Recalling the definition of $\cG$, we have the following observation,
\begin{align}
\label{eq:ADMM-cg} 
    \cG(\frak{a},\frak{b}) &= - \inf_{m \geq 0} \int_{\dom_{T}} m(t,x) \Bigl(\frak{a}(t,x) - H\bigl(x, m(t,x), \frak{b}(t,x) \bigr) \Bigr) \rd x \, \rd t \nonumber \\
    & = - \int_{\dom_{T}} \inf_{m \geq 0} \Bigl [m \bigl(\frak{a}(t,x) - H\bigl(x, m, \frak{b}(t,x) \bigr) \bigr) \Bigr] \rd x \, \rd t \\
    & = -\int_{\dom_{T}} \cK\bigl(\frak{a}(t,x), \frak{b}(t,x)\bigr) \rd x \, \rd t,
\end{align}
where we denoted $\cK\bigl(a, b\bigr) = \inf_{m \geq 0} \Bigl [m \bigl(\frak{a} - H\bigl(x, m, \frak{b} \bigr) \bigr) \Bigr]$.

Therefore, the form of $\cG$ depends on $\cK$. In general, we do not know any closed form of $\cK$ in terms of $\frak{a}$ and $\frak{b}$. However, for several Hamiltonian functions $H(x,m,p)$ of interest, we can derive such closed-form solution. %
Here we demonstrate some of the calculations to deliver a general idea. 

\begin{example}[Mean-field Aversion]
Consider the running cost $f(m,v) = \frac{1}{4}\|v\|^{2} + m$ and drift $b(x,m,v) = v$. The corresponding Hamiltonian is $H(x,m,p)$ = $\|p\|^{2} - m$. Then
\begin{align}
\label{eq:k-aversion}
    \cK(\frak{a},\frak{b}) 
    = & \inf_{m \geq 0} \; m^{2} + (\frak{a} - \|\frak{b}\|^{2}) m %
    = 
    \begin{cases}
     -\frac{1}{4}\bigl(\frak{a} - \|\frak{b}\|^{2} \bigr)^{2} \quad  &\hbox{ if } \frak{a}-\|\frak{b}\|^{2} \geq 0 \\
        0  \quad &\hbox{ otherwise. } 
    \end{cases}
\end{align}
\end{example}

\begin{example}[Mean-field Maximum Entropy]
Consider the running cost $f(m,v) = \frac{1}{4}\|v\|^{2} + \log(m)$ and drift $b(x,m,v) = v$, the corresponding Hamiltonian is $H(x,m,p) = \frac{1}{2} \|p\|^{2} - \log(m)$, then
\begin{align}
\label{eq:k-max-entropy}
    \cK(\frak{a},\frak{b}) 
    = & \inf_{m \geq 0} \; m\log m + (\frak{a} - \frac{1}{2}\|\frak{b}\|^{2}) m %
    = - \exp \bigl(\frac{1}{2}\|\frak{b}\|^{2} - \frak{a} - 1\bigr).
\end{align}
\end{example}

\begin{example}[Continuous Optimal Transport]
Consider the running cost $f(m,v) = \frac{1}{2}\|v\|^{2}$ and drift $b(x,m,v) = v$, the corresponding Hamiltonian is $H(x,m,p) = \frac{1}{2}\|p\|^{2}$, then
\begin{align}
\label{eq:k-ot}
\cK(\frak{a},\frak{b}) = & \inf_{m \geq 0} \; m(\frak{a} - \frac{1}{2}\|\frak{b}\|^{2}) %
= 
\begin{cases}
   0 \quad \quad &\hbox{ if } \frak{a}-\frac{1}{2}\|\frak{b}\|^{2} \geq 0 \\
   -\infty \quad &\hbox{ otherwise.}
\end{cases}
\end{align}
\end{example}

\begin{example}[Mean-field Congestion]
Consider the running cost $f(m,v) = \frac{1}{4}m\|v\|^{2}$ and drift $b(x,m,v) = v$, the corresponding Hamiltonian is $H(x,m,p) = \frac{\|p\|^{2}}{m}$, then
\begin{align}
\label{eq:k-congestion}
    \cK(\frak{a},\frak{b}) = & \inf_{m \geq 0} \; m(\frak{a} - \frac{\|\frak{b}\|^{2}}{m}) 
= 
\begin{cases}
   -\|\frak{b}\|^{2} \quad \quad &\hbox{ if }  \frak{a} \geq 0 \\
   -\infty \quad &\hbox{ otherwise. }
\end{cases}
\end{align}
\end{example}

From these examples, it can be seen that $\cK$ has a closed form for many smooth Hamiltonian. However, the existence of closed form expressions of $\cK$ alone is not enough for making the training loss tractable. In the example of \eqref{eq:k-ot} and \eqref{eq:k-congestion}, $\cK$ takes values $-\infty$, which makes $\cG$ singular and computationally intractable at some points. Moreover, gradient-based training cannot be carried out successfully in the presence of infinite values as well. The presence of $-\infty$ in $\cK$ is due to the degeneracy of $H$ in terms of order in $m$. 
For cases with singular $\cK$ and $\cG$, we need an additional trick to tackle this issue.

Recall that in Algorithm \ref{algo:vanilla-admm}, the update for function $q$ is given by,
\begin{align}
\label{eq:admm-minmax-1}
    q^{(k)} 
    = & \underset{q: \dom_{T} \rightarrow \RR^{d+1}}{\mathop{\arg\min}} \; \mathcal{G}(q) + \inp{\lambda^{(k-1)}}{q} + \frac{r}{2} \bigl\| \Lambda u^{(k)} - q\bigr\|^{2} \nonumber 
    \\
    = & \underset{q: \dom_{T} \rightarrow \RR^{d+1}}{\mathop{\arg\min}} \; \int_{\dom_{T}} \left(-\cK(q(t,x)) +  \inp{\lambda^{(k-1)}(t,x)}{q(t,x)} + \frac{r}{2} \bigl\| \Lambda u^{(k)}(t,x) - q(t,x) \bigr\|^{2} \right)\; \rd x \; \rd t
\end{align}
Since we don't have any additional constraint on the value of $q$, the function $q^{(k)}$ that minimizes the integral in \eqref{eq:admm-minmax-1} should minimize the integrand point wisely. Therefore, it holds that:
\begin{align}
\label{eq:admm-minmax-2}
    q^{(k)}(t,x) = \; & \underset{q \in \RR^{d+1}}{\mathop{\arg\min}} \quad \left( -\cK(q) + \inp{\lambda^{(k-1)}(t,x)}{q} + \frac{r}{2} \bigl\| \Lambda u^{(k)}(t,x) - q \bigr\|^{2} \right) \nonumber 
    \\
    = \; & \underset{q \in \RR^{d+1}}{\mathop{\arg\min}} \sup_{m \geq 0} \; m\bigl(H(x, m, q_{2}) - q_{1}\bigr) + \inp{\lambda^{(k-1)}(t,x)}{q} + \frac{r}{2} \bigl\| \Lambda u^{(k)}(t,x) - q \bigr\|^{2},
\end{align}
where we denote $q = (q_1, q_2)$, with $q_{1} \in \RR$ and $q_2 \in \RR^{d}$. For fixed $\lambda^{(k-1)}$ and $u^{(k)}$, let us define $\cL(q, m) = m\bigl(H(x, m, q_{2}) - q_{1}\bigr) + \inp{\lambda^{(k-1)}(t,x)}{q} + \frac{r}{2} \bigl\| \Lambda u^{(k)}(t,x) - q \bigr\|^{2}$. Under certain conditions, the following minimax equality holds for the right-hand side of~\eqref{eq:admm-minmax-2},
\begin{align}
\label{eq:admm-minmax-3}
\inf_{q \in \RR^{d+1}} \sup_{m \geq 0} \cL(q, m) = \sup_{m \geq 0} \inf_{q \in \RR^{d+1}} \cL(q, m)
\end{align}
Therefore, we can exploit \eqref{eq:admm-minmax-3} to address the issues of singular $\cK$. We notice that for the cases of congestion and general linear quadratic, $\cL(q,m)$ is quadratic in $q$ for fixed $m$. This means that we can solve explicitly in a closed form for $\inf_{q} \cL(q,m)$ for fixed $m$, then we solve for the maximization problem over $m \geq 0$. In this way, we can effectively avoid the issues generated by infinity values in $\cK$. \\

Moreover, using this trick, we directly obtain the value for $q^{(k)}$ at any $(t,x)$ based on the value of $u^{(k)}(t,x)$ and  $\lambda^{(k-1)}(t,x)$. Therefore, we can skip the neural network training in Algorithm \ref{algo:DeepADMM} for function $q$ in each DeepADMM iteration. Instead, we compute the values for $q^{(k)}(t,x)$ directly following the above procedure whenever the values are needed to compute $\cL^{(u)}$ and $\cL^{(\lambda)}$.

\section{Some loss plots}

Contrary to the LQ case, we do not have any benchmark solution to compare our numerical results to for the congestion case. We already explained why the methods gave consistent results, namely that we see the density spreading before reforming the terminal distribution, in contrast with the LQ case where the shape of the population remains the same along the trajectory (see~\ref{subsec:cong_num}). In this section, we provide some plots showing the evolution of the different losses through the training of the neural networks. 

\begin{figure}[t]
    \centering
    \includegraphics[width=0.5\textwidth]{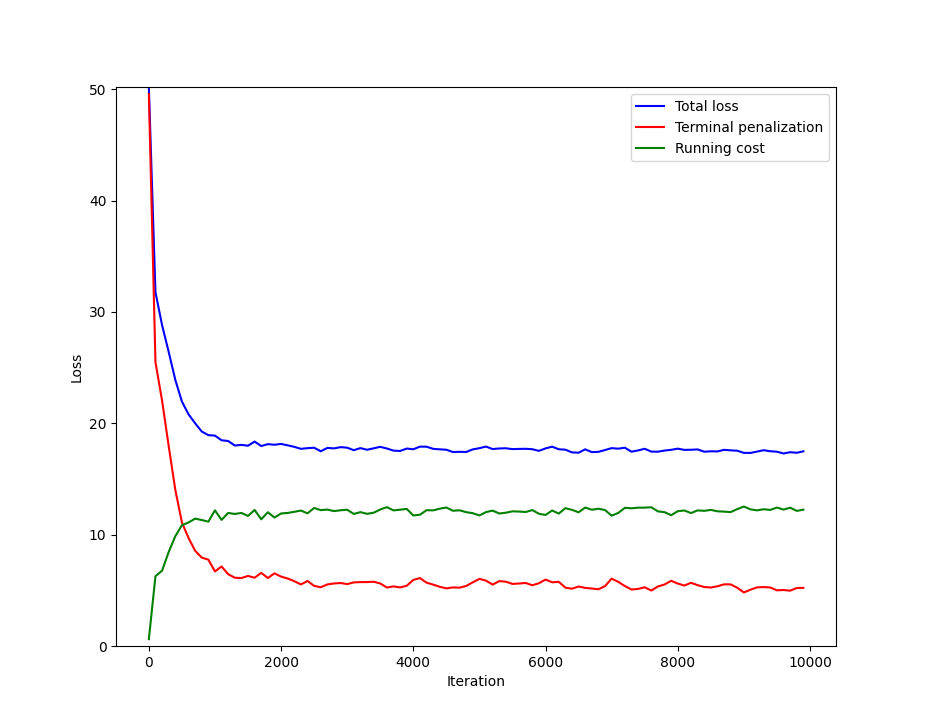}
    \caption{Evolution of the loss in the congestion case in dimension 5 for method 1. }
    \label{fig:cong5D_loss_m1}
\end{figure}

Figure~\ref{fig:cong5D_loss_m1} shows the losses for Method~1. The total loss is the sum of the ``running cost'' loss and the ``terminal penalization'' loss multiplied by $C_W$. In this test case, we took $C_W=10$, meaning that the penalization is ten times the Wasserstein distance between the effective terminal distribution of agents and the desired one. We observe that for the first iterations, the distance is high and the running cost almost zero, which comes from the fact that the first try is not to move. Then the algorithm makes the penalization decrease by paying a trade-off in the form of the running cost, and seems to reach a plateau.

\begin{figure}[t]
    \centering
    \includegraphics[width=0.5\textwidth]{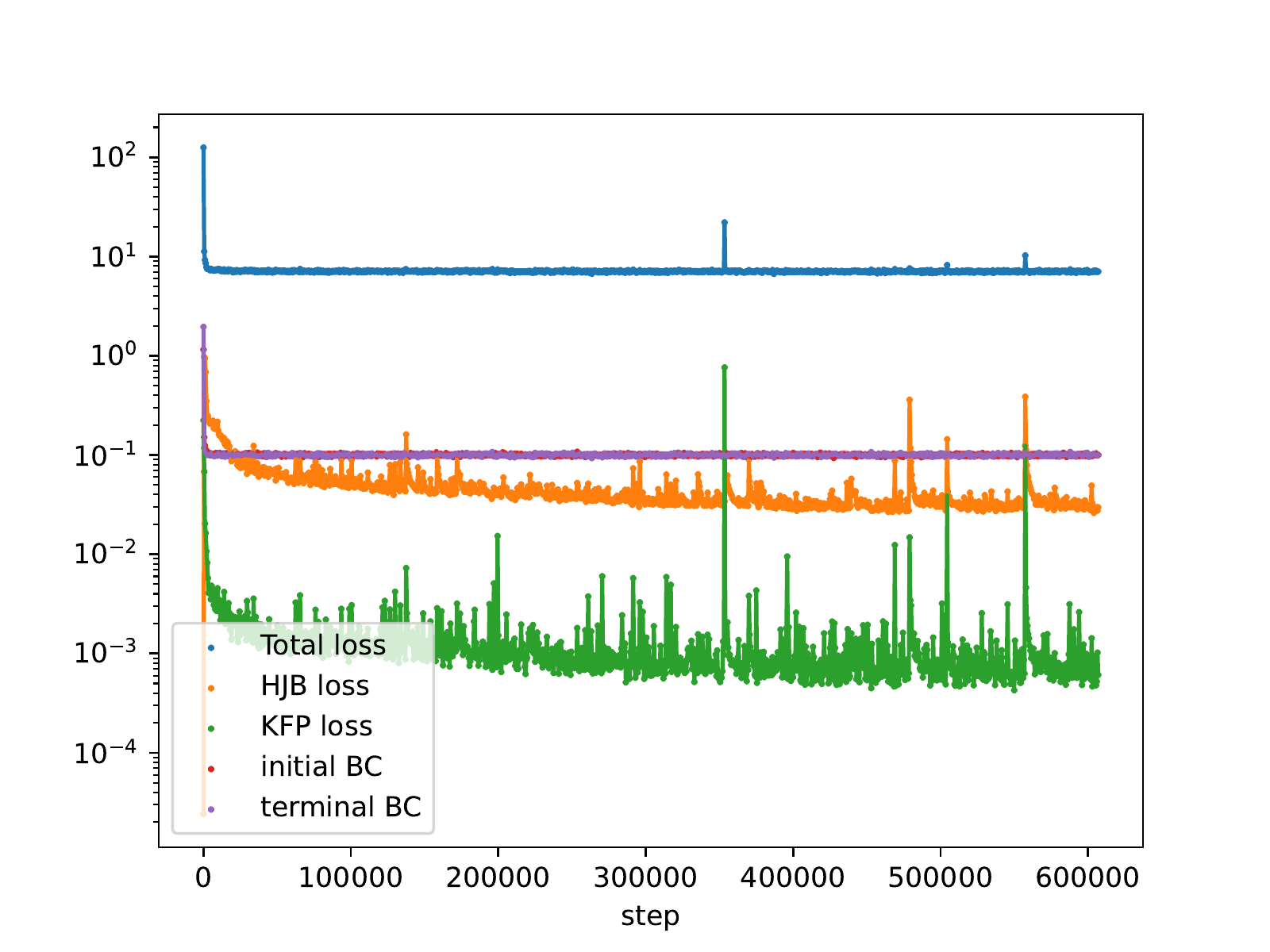}
    \caption{Evolution of the loss in the congestion case in dimension 5 for method 2. }
    \label{fig:lq_test1_evol_M2}
\end{figure}

Figure~\ref{fig:lq_test1_evol_M2} shows the losses for Method~2. The ``HJB loss'' is the squared $L^2$ residual for the HJB equation. The ``KFP loss'' is the squared $L^2$ residual for the KFP equation. The ``initial BC'' loss and the ``terminal BC'' loss respectively correspond to the squared $L^2$ error on the initial and terminal distributions. 
The ``total loss'' is the sum of the other losses, up to multiplicative weights. 

\begin{figure}[t]
    \centering
    \includegraphics[width=0.5\textwidth]{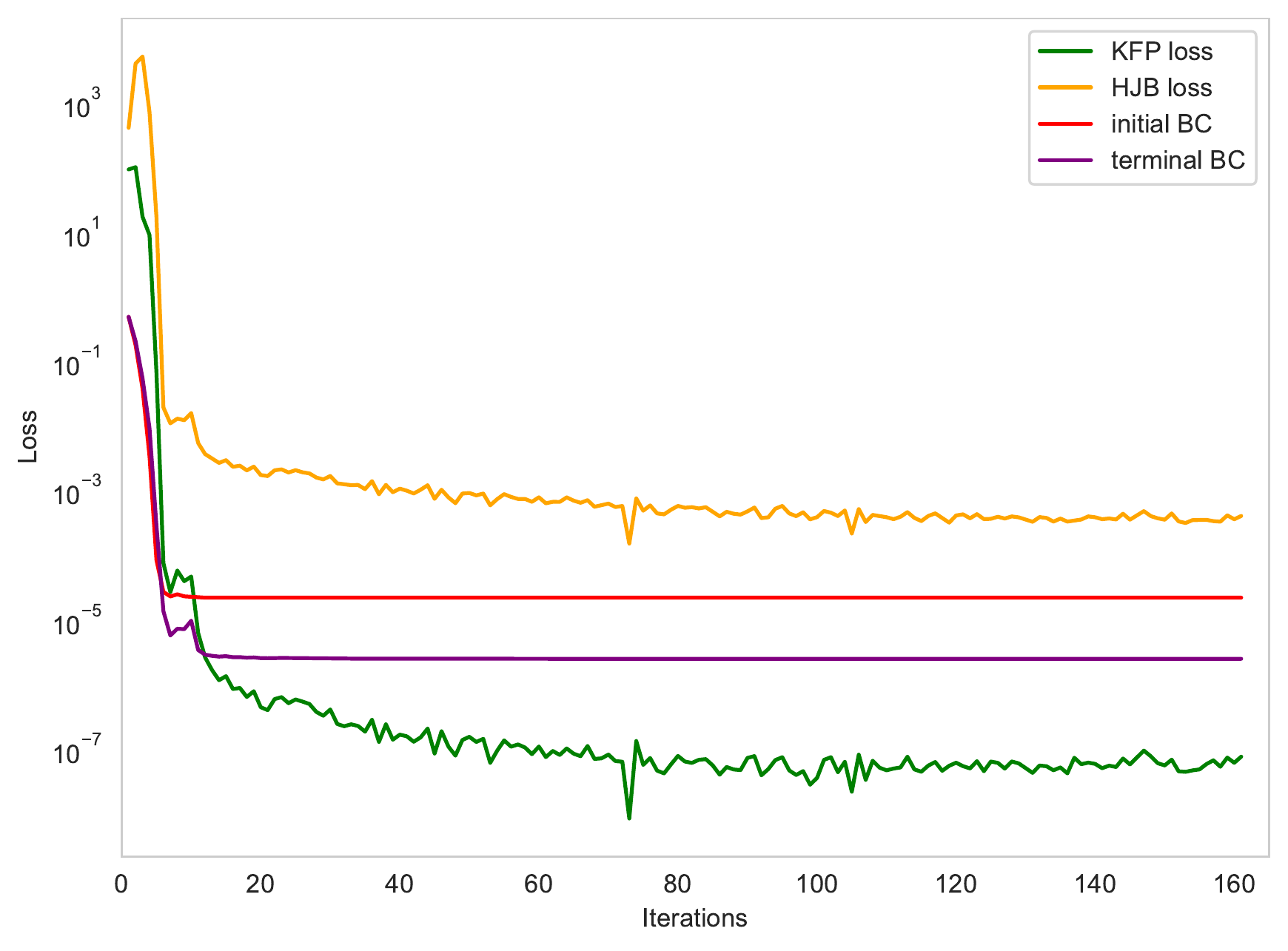}
    \caption{Evolution of the loss in the congestion case in dimension 5 for method 3. }
    \label{fig:cong5D_loss_m3}
\end{figure}

Figure~\ref{fig:cong5D_loss_m3} shows, for Method~3, the squared $L^2$ residuals for the HJB and KFP equations, as well as the squared $L^2$ loss for the initial and terminal conditions. Note that these losses are not directly minimized during the algorithm of Method~3, but they are minimized as a by-product of the iterations. 

\end{document}